\documentclass[english]{article}
\usepackage{lmodern}
\usepackage{lmodern}
\usepackage[T1]{fontenc}
\usepackage[utf8]{inputenc}
\usepackage[a4paper]{geometry}
\usepackage[active]{srcltx}
\usepackage{babel}
\usepackage{mathrsfs}
\usepackage{mathtools}
\usepackage{tcolorbox}
\usepackage{amsmath}
\usepackage{amsthm}
\usepackage{amssymb}
\usepackage{microtype}
\usepackage[bookmarks=false,
 breaklinks=false,pdfborder={0 0 1},backref=section,colorlinks=false]
 {hyperref}
\hypersetup{hidelinks} 

\makeatletter

\usepackage{xcolor}
\usepackage{lineno}

\usepackage{mathrsfs}
\usepackage{tikz-cd}

\usepackage{amsthm}
\theoremstyle{plain}
\newtheorem{theorem}{Theorem}\newtheorem{lemma}[theorem]{Lemma}\newtheorem{proposition}[theorem]{Proposition}\newtheorem{corollary}[theorem]{Corollary}\theoremstyle{definition}
\newtheorem{definition}[theorem]{Definition}\theoremstyle{remark}
\newtheorem{remark}[theorem]{Remark}\theoremstyle{definition}
\newtheorem{example}[theorem]{Example}

\usepackage{wasysym}
\usepackage{graphicx}
\usepackage{upgreek}
\usepackage{hyphenat}

\usepackage{dsfont}

\usepackage{babel}

\makeatother

\begin{document}
	
%\linenumbers	
	
\title{Monotone extensions into compact target sets in dual systems}
\author{M.D. Voisei\thanks{Corresponding author.
		Email: \texttt{mdvoisei@yahoo.com}}\\
	Baltimore, Maryland, U.S.A.
	}
\date{{}}

\maketitle
 
\begin{abstract}
	We study monotone extension problems in dual systems without assuming
	separation. Our compact target-set formulation includes multivalued
	operators as a special case and allows the initial set to be
	nonmonotone. Using methods from nonlinear and convex analysis, we
	establish existence results and characterize when a compact target
	contains a point monotonically related to the initial set. We apply
	these results to variational inequalities and give an example in which
	an auxiliary topology yields a sharp threshold for existence. We also
	establish a compact-target criterion for global minimization of
	differentiable functions without assuming convexity. Finally, for a
	Cauchy location objective with a variable reference point, we prove
	the existence of an estimate that uniquely minimizes the objective
	when the reference point is fixed at the value determined by that
	estimate.
\end{abstract}

\strut


\section{Introduction}

We work in a dual system $(X,Y,\langle\cdot,\cdot\rangle)$, where
$X$ and $Y$ are real vector spaces paired by a bilinear form $\langle\cdot,\cdot\rangle:X\times Y\to\mathbb{R}$.
Given subsets $G,S\subset X\times Y$, we study the problem of extending
$G$ monotonically into the target set $S$. More precisely, for $A\subset X\times Y$,
define $A^{+}$ as the set of all points that are monotonically related
to $A$: 
\[
A^{+}:=\left\{ (x,y)\in X\times Y\;\mid\forall(u,v)\in A,\;\langle x-u,y-v\rangle\geq0\right\} .
\]

The main problem considered in this paper is to identify conditions
under which 
\[
G^{+}\cap S\neq\emptyset,
\]
or equivalently, to find points in $S$ that are monotonically related
to $G$. Unlike the classical graph-over-compact-domain formulation,
the compactness in our approach is imposed on the target set $S$,
not necessarily on the domain of $G$ or on the graph of an operator.

When $T:X\rightrightarrows Y$ is a multivalued operator, the preceding
formulation applies by taking 
\[
G=\operatorname{Graph}T:=\{(x,y)\in X\times Y\mid y\in T(x)\}.
\]

Thus, for a prescribed target set $S\subset X\times Y$, the extension
problem asks whether there exists $(x,y)\in S$ that is monotonically
related to the graph of $T$. The target set may therefore encode
additional restrictions on the desired extension point, such as compactness,
convexity, dissipativity, range constraints, variational-inequality
constraints, or fixed-point constraints.

This formulation is motivated by the observation that, in a monotone extension problem, the set against which monotonicity is tested and the set in which the extension point is sought play different roles and need not be subject to the same assumptions. In the compactness-based existence results considered here, compactness can be imposed on the latter rather than on the domain of the initial set or operator. Separating these two roles leads to the problem \(G^{+}\cap S\neq\varnothing\), in which the comparison set \(G\) may be noncompact and even nonmonotone, while the target set \(S\) carries the compactness and geometric assumptions used to obtain existence. This formulation also allows constraints on the desired extension point to be incorporated directly into the choice of \(S\).

The problem of monotone extension has been addressed in only a few
results in the literature, beginning with the classical work of Debrunner-Flor
\cite{MR0170189} and Browder \cite[Theorem 8]{MR0229101} and continuing
in several related later results. These results typically formulate
the problem for an operator whose domain is contained in a compact
convex set, under continuity, semicontinuity, or closedness assumptions
imposed through a suitable topology; see, for example, \cite[Theorem 12]{MR190689},
\cite[Theorem 2.7]{MR721236}, \cite[Theorem 3]{MR1938858}.

The contribution of the paper is twofold. First, we replace the classical
graph-over-compact-domain setting by a target-set formulation: one
seeks a point in a compact subset of $X\times Y$ that is monotonically
related to a prescribed set $G$. In this formulation, $G$ need not
be monotone and need not have compact domain. This distinction is
essential in the main minimax theorem, where $S$ is compact, convex,
and dissipative, while $G$ enters through the convexification of
its quadratic form.

Second, the results are obtained in general dual systems, which are not assumed to be Hausdorff separated, rather than only in the standard setting of locally convex spaces and their continuous duals. The motivation for allowing non-separated dual systems is that the monotone-extension existence problem depends only on the effective duality induced by the bilinear pairing, and therefore does not intrinsically require a separation assumption. Indeed, the nontrivial annihilator directions are invisible to the pairing and can be removed by passing to the corresponding separated quotient without changing the existence problem. Accordingly, the non-separated case is handled by passing to appropriate separated quotients when needed, while the statements remain intrinsic to the original dual system (see Section 4).

Our proofs combine classical techniques with additional convex-analytic
tools. As in the existing literature, some arguments rely on fixed-point
theorems, such as Brouwer's theorem and the Ky Fan-Glicksberg theorem,
and on partitions of unity associated with finite open covers of compact
spaces. The new ingredients are minimax principles and Fenchel-Rockafellar
duality, which make it possible to treat the extension problem in
the compact target-set formulation described above. In particular,
the central convex-analytic criterion expresses the nonemptiness of
$G^{+}\cap S$ through inequalities involving the Fitzpatrick-type
functions associated with $G$, thereby connecting monotone extension
to convex duality. 

The resulting theory yields both existence theorems and convex-analytic
criteria for monotone extension into compact target sets. Several
results give necessary and sufficient conditions for the nonemptiness
of $G^{+}\cap S$, while others provide directly applicable sufficient
conditions; see for example Theorems \ref{E}, \ref{H}, \ref{IM},
\ref{BRS}, and \ref{DFE} and their consequences below. 

Theorems \ref{E} and \ref{H} provide the minimax and compact-target
extension principles; Theorem \ref{IM} gives an intersection criterion
under upper semicontinuity assumptions; Theorem \ref{BRS} recovers
and extends Browder-type extension results; and Theorem \ref{DFE}
gives a Debrunner-{}-Flor-type fixed-point consequence.

These extension principles also lead to existence results for variational
inequalities; see, for example, Corollary \ref{VIS} and Theorems
\ref{VI0}, \ref{VID} below. 

The compact-target formulation is particularly useful for these variational-inequality applications because it allows the set in which the solution is sought to be distinct from the set of comparison points. In Corollary \ref{VIS} this distinction is represented by the target set \(S\) and the comparison set \(C\), respectively, while in Theorems \ref{VI0} and \ref{VID} the solution is sought in \(K\) and the comparison points belong to \(\operatorname{conv} C\). The classical situation in which the solution and comparison sets coincide is recovered as a particular case.

\section{Preliminaries}

A dual system $(X,Y,\langle\cdot,\cdot\rangle)$ is said to be \emph{separated
in} $X$ if, for every $x\in X$ with $x\neq0$, there exists $y\in Y$
such that $\langle x,y\rangle\neq0$. This condition is equivalent
to $(X,\sigma(X,Y))$ being Hausdorff. Here $\sigma(X,Y)$ denotes
the weakest topology on $X$ for which all linear forms $\{X\ni x\to\langle x,y\rangle\in\mathbb{R}\mid y\in Y\}$
are continuous. This topology is locally convex and called \emph{the
weak topology of} $X$ with respect to the duality $(X,Y,\langle\cdot,\cdot\rangle)$.

Separation in $Y$ is defined in a similar manner. A dual system $(X,Y,\langle\cdot,\cdot\rangle)$
is called \emph{separated} if it is separated both in $X$ and $Y$.

To a multi-valued operator (or multi-function) $T:X\rightrightarrows Y$
we associate its \emph{values}: $T(x)\subset Y$, $x\in X$, \emph{inverse:}
$T^{-1}:Y\rightrightarrows X$, $\operatorname*{Graph}(T^{-1})=\{(y,x)\mid(x,y)\in\operatorname*{Graph}T\}$,
\emph{domain}: $D(T):=\{x\in X\mid T(x)\neq\emptyset\}=\Pr\nolimits_{X}(\operatorname*{Graph}T)$,
\emph{range}: $R(T):=\{y\in Y\mid y\in T(x)\ \mathrm{for\ some}\ x\in X\}=\Pr\nolimits_{Y}(\operatorname*{Graph}T)$,
\emph{direct--image:} $T(A)=\cup\{T(x)\mid x\in A\}\subset Y$, for
$A\subset X$, \emph{inverse--image:} $T^{-1}(B)=\{x\in X\mid T(x)\cap B\neq\emptyset\}\subset X$,
for $B\subset Y$. Here $\Pr_{X}(x,y):=x$ and $\Pr_{Y}(x,y):=y$,
$(x,y)\in X\times Y$, are the projections of $X\times Y$ onto $X$
and $Y$, respectively.

A set-valued mapping $T:X\rightrightarrows Y$ is \emph{monotone}
if 
\[
\forall(x_{1},y_{1}),(x_{2},y_{2})\in\operatorname*{Graph}T,\ \left\langle x_{1}-x_{2},y_{1}-y_{2}\right\rangle \geq0.
\]
A subset $A\subset X\times Y$ is \emph{monotone} if it is the graph
of a monotone mapping, i.e., $A=\operatorname*{Graph}T$, for some
monotone $T:X\rightrightarrows Y$, or equivalently, for every $(x_{1},y_{1}),(x_{2},y_{2})\in A$,
$\left\langle x_{1}-x_{2},y_{1}-y_{2}\right\rangle \geq0$.

A set $A\subset X\times Y$ is \emph{dissipative} if, for every $(x_{1},y_{1}),(x_{2},y_{2})\in A$,
$\left\langle x_{1}-x_{2},y_{1}-y_{2}\right\rangle \le0$, or equivalently,
the set 
\begin{equation}
\neg A:=\{(x,-y)\mid(x,y)\in A\}\label{eq:-38}
\end{equation}
is monotone.

Given a dual system $(X,Y,c=\langle\cdot,\cdot\rangle)$, the product
space $Z:=X\times Y$ forms a \emph{natural dual system} $(Z,Z,\cdot)$
where, for $z=(x,y),\ w=(u,v)\in Z$, the bilinear form $\cdot:Z\times Z\to\mathbb{R}$
is defined as 
\begin{equation}
z\cdot w:=\langle x,v\rangle+\langle u,y\rangle.\label{eq:-4}
\end{equation}

Throughout the paper, $Z$ is endowed with a locally convex topology
$\mu$ compatible with this duality, in the sense that every $\mu-$continuous
linear functional on $Z$ is represented by an element of $Z$ through
the pairing ``$\cdot$''. Equivalently 
\[
\ell\in(Z,\mu)^{*}\ \Rightarrow\ \ell(z)=z\cdot w,\ \mathrm{for\ some}\ w\in Z\ \mathrm{and\ all}\ z\in Z.
\]

One example of such a topology is $\mu=\sigma(Z,Z)=\sigma(X,Y)\times\sigma(Y,X)$.

With respect to the natural dual system $(Z,Z,\cdot)$, the conjugate
of a function $f:Z\rightarrow\overline{\mathbb{R}}$ is denoted by
$f^{\square}:Z\rightarrow\overline{\mathbb{R}}$ and is given by 
\begin{equation}
f^{\square}(z)=\sup\{z\cdot z^{\prime}-f(z^{\prime})\mid z^{\prime}\in Z\}.\label{eq:-7}
\end{equation}
By the biconjugate formula, $f^{\square\square}=\operatorname*{cl}_{(\mu)}\operatorname*{conv}f$
whenever $f^{\square}$ (or equivalently, $\operatorname*{cl}\operatorname*{conv}f$)
is proper. Here $\operatorname*{cl}_{\mu}f$ denotes the $\mu-$\emph{lower
semicontinuous hull} of $f$, which is the greatest $\mu-$lower semicontinuous
function not exceeding $f$; $\operatorname*{conv}f$ is the \emph{convex
hull} of $f$, i.e., the greatest convex function not exceeding $f$.

For $A\subset Z$ we associate the following functions: $\iota_{A},c_{A},\varphi_{A},\psi_{A}:Z\rightarrow\overline{\mathbb{R}}$
defined by 
\begin{itemize}
\item $\iota_{A}(z)=0$, if $z\in A$, $\iota_{A}(z)=+\infty$, otherwise; 
\item $c_{A}:=c+\iota_{A}$; 
\item $\varphi_{A}:=c_{A}^{\square}$ -- called the \emph{Fitzpatrick function}
of $A$; 
\item $\psi_{A}:Z\rightarrow\overline{\mathbb{R}}$, $\psi_{A}:=\varphi_{A}^{\square}=c_{A}^{\square\square}$. 
\end{itemize}
Similarly, for a multi\-function $T:X\rightrightarrows Y$ $c_{T}:=c_{\operatorname*{Graph}T}$,
$\psi_{T}:=\psi_{\operatorname*{Graph}T}$, and $\varphi_{T}:=\varphi_{\operatorname*{Graph}T}$
is the \emph{Fitzpatrick function} of $T$.

For $T:X\rightrightarrows Y$ we define $T^{+}:X\rightrightarrows Y$
via $\operatorname*{Graph}(T^{+}):=(\operatorname*{Graph}T)^{+}$,
or equivalently 
\begin{equation}
\operatorname*{Graph}(T^{+})=[\varphi_{T}\le c]:=\{z\in Z\mid\varphi_{T}(z)\le c(z)\}.\label{eq:-3}
\end{equation}

In general, for a given topological vector space $(E,\mu)$, a subset
$A\subset E$, and functions $f,g:E\rightarrow\overline{\mathbb{R}}$
we introduce the following notations: 
\begin{itemize}
\item $\operatorname*{conv}A$ -- the \emph{convex hull} of $A$; $x\in\operatorname*{conv}A$
if and only if $(x,1)={\displaystyle \sum_{k=1}^{n}}\lambda_{k}(x_{k},1)$
for some $n\in\mathbb{N},$ $n\ge1$, $\{\lambda_{k}\}_{k\in\overline{1,n}}\subset\mathbb{R}^{+}:=[0,+\infty)$
and $\{x_{k}\}_{k\in\overline{1,n}}\subset A$. 
\item $\operatorname*{aff}A$ -- the \emph{affine hull} of $A$. A point
$x\in\operatorname*{aff}A$ if and only if $(x,1)={\displaystyle \sum_{k=1}^{n}}\lambda_{k}(x_{k},1)$
for some $n\in\mathbb{N},$ $n\ge1$, $\{\lambda_{k}\}_{k\in\overline{1,n}}\subset\mathbb{R}$
and $\{x_{k}\}_{k\in\overline{1,n}}\subset A$. 
\item $\operatorname*{lin}A=\operatorname*{aff}(A\cup\{0\})$ -- the \emph{linear
hull} of $A$, which is the smallest subspace of $E$ containing $A$. 
\item $[f\le g]$ $:=\{x\in E\mid f(x)\leq g(x)\}$. The sets $[f=g]$,
$[f<g]$, $[f>g]$, $[f\ge g]$ are defined in a similar manner. 
\item $\operatorname*{cl}_{\mu}A$ -- the $\mu-$\emph{closure} of $A$.
A point $x\in\operatorname*{cl}_{\mu}A$ if and only if $x_{i}\stackrel{\mu}{\to}x$,
for some net $(x_{i})_{i}\subset A$. Here ``$\stackrel{\mu}{\to}$''
denotes the convergence of nets in the $\mu$ topology. 
\end{itemize}
We omit the explicit $\mu-$notation when the topology is implicitly
understood.

\medskip{}

For properties of monotone and representable operators we refer to
\cite{MR2453098,MR2594359,MR2577332}. Some of these frequently used
properties are as follows: 
\begin{itemize}
\item $T$ is monotone if and only if $\operatorname*{conv}c_{T}\ge c$
if and only if $\psi_{T}\ge c$. In this case $\operatorname*{Graph}T\subset[\psi_{T}=c]$. 
\item $c_{T}$ is convex if and only if $T$ is monotone and $\operatorname*{Graph}T$
is convex. 
\end{itemize}
\medskip{}

In this paper, the following conventions are adopted: $\inf\emptyset=+\infty$,
$\sup\emptyset=-\infty$, $\infty-\infty=\infty$, $0\cdot\pm\infty=0$.

\eject

\section{Analysis and interpretation of the Debrunner--Flor and Browder results}

Recall the following results of Debrunner-Flor and Browder:

\begin{theorem} \label{DF} (Debrunner-Flor \cite{MR0170189}) Let
$X$ and $Y$ be real topological vector spaces, where $X$ is a Hausdorff
locally convex space. Let $M\subset X\times Y$ be monotone with respect
to a continuous bilinear form $\langle\cdot,\cdot\rangle$ on $X\times Y$.
Suppose that the domain of $M$ is contained in a compact convex subset
$A$ of $X$. If $\varphi:A\to Y$ is continuous then, there exists
$x\in A$ such that the extended set $M\cup\{(x,\varphi(x))\}$ remains
monotone. \end{theorem}

\begin{theorem} \label{BR} (Browder \cite[Theorem 8]{MR0229101})
Let $K$ be a compact convex subset of the real Hausdorff topological
vector space $(E,\tau)$, $(F,\mu)$ a topological vector space, with
a bilinear pairing given between $E$ and $F$ to the reals which
we denote by $\langle w,u\rangle$ for $w$ in $F$ and $u$ in $E$.
We suppose that the mapping of $K\times F$ into $\mathbb{R}$ which
carries $(u,w)$ into $\langle w,u\rangle$ is continuous. Let $f$
be a continuous mapping of $K$ into $F$ and let $G$ be a monotone
subset of $K\times F$. Then there exists an element $u_{0}$ of $K$
such that $(u_{0},f(u_{0}))$ is monotonically related to $G$, or
equivalently, $G\cup\{(u_{0},f(u_{0}))\}$ is monotone. \end{theorem}

The conclusions of Theorems \ref{DF}, \ref{BR} share a similar structure;
in their respective notations, 
\begin{equation}
M^{+}\cap\operatorname*{Graph}\varphi\neq\emptyset,\ G^{+}\cap\operatorname*{Graph}f\neq\emptyset.\label{eq:-2}
\end{equation}

In the proof of Debrunner-Flor Theorem, it suffices that the restricted
bilinear form $\langle\cdot,\cdot\rangle$ be continuous on $A\times Y$,
where $A$ carries the topology inherited from $X$. Therefore, Theorem
\ref{DF} can be viewed as a special case of Theorem \ref{BR}, whose
formulation has the additional advantage that $(E,\tau)$ need only
be a Hausdorff topological vector space rather than a locally convex
one.

If the function $\varphi$ in Theorem \ref{DF} is constant, then
the proof only requires that, for every $y\in Y$, $\langle\cdot,y\rangle$
be continuous on $A$. These observations motivate the following notion.

\begin{definition} \label{w} Given a dual system $(X,Y,c)$ and
$\emptyset\neq B\subset Y$, we define $\sigma(X,B)$ to be the locally
convex topology on $X$ generated by the family of seminorms $\{|c(\cdot,y)|\mid y\in B\}$.
This is the weakest topology on $X$ for which all linear forms $\{X\ni x\to c(x,y)\mid y\in B\}$
are continuous. For $\emptyset\neq A\subset X$ we define $\sigma(Y,A)$
similarly. \end{definition}

For $\emptyset\neq G\subset Z$ we define $\sigma(Z,G)$ analogously
to Definition \ref{w}, relative to the dual system $(Z,Z,\cdot)$
defined in (\ref{eq:-4}).

\begin{remark} \label{RI} Under the assumptions and notations of
Theorem \ref{BR} 
\begin{itemize}
\item $K$ is $\sigma(E,F)-$compact. Indeed, due to the $\tau\times\mu-$continuity
of $\langle\cdot,\cdot\rangle:K\times F\to\mathbb{R}$, the topology
$\sigma(E,F)$ is weaker than $\tau$ on $K$. Since $K$ is $\tau-$compact,
it is also compact with respect to any topology weaker than $\tau$,
including $\sigma(E,F)$. 
\item $\operatorname*{Graph}f$ is $\tau\times\mu-$compact and the topology
$\tau\times\mu$ is Hausdorff separated on $\operatorname*{Graph}f$.
Indeed, concerning the $\tau\times\mu-$Hausdorff separation of $\operatorname*{Graph}f$,
it suffices to note that if $(x,f(x))\neq(y,f(y))$, then necessarily
$x\neq y$, and we use the Hausdorff separation of $(E,\tau)$. It
is easily checked that, due to the continuity of $f$ and $\tau-$compactness
of $K$, $\operatorname*{Graph}f$ is $\tau\times\mu-$compact which
allows us to conclude that $\operatorname*{Graph}f$ is compact with
respect to any topology which is weaker than $\tau\times\mu$ on $\operatorname*{Graph}f$. 
\item For example, on $K\times F$, the product topology $\sigma(E,F)\times\sigma(F,K)$
is weaker than $\tau\times\mu$. This follows from the fact that $\sigma(E,F)$
is weaker than $\tau$ on $K$, and $\sigma(F,K)$ is weaker than
$\mu$ on $F$, both due to the $\tau\times\mu-$continuity of $\langle\cdot,\cdot\rangle:K\times F\to\mathbb{R}$.
In general, $\sigma(E\times F,G)$ is weaker than $\sigma(E,F)\times\sigma(F,K)$
on $E\times F$ since $G\subset K\times F$, where $(E\times F,E\times F,\cdot)$
is the naturally associated dual system, defined as in (\ref{eq:-4}).
In conclusion, on $K\times F$, and consequently on $\operatorname*{Graph}f\subset K\times F$,
the topology $\tau\times\mu$ is stronger than $\sigma(E,F)\times\sigma(F,K)$,
which in turn is stronger than $\sigma(E\times F,G)$. Therefore $\operatorname*{Graph}f$
is compact with respect to both $\sigma(E,F)\times\sigma(F,K)$ and
$\sigma(E\times F,G)$. 
\end{itemize}
\end{remark}

\section{The separated version of a dual system }

Let $N$ be a vector subspace of a vector space $X$ and let $X/N$
be the \emph{quotient space of $X$ modulo} $N$ consisting of equivalence
classes $\pi_{N}(x):=\hat{x}:=x+N$, $x\in X$, where $\pi_{N}:X\to X/N$
is the canonical \emph{quotient map} of $X$ onto $X/N$. If $\tau$
is a linear topology on $X$, then \emph{quotient topology} $\tau_{N}$
on $X/N$ is defined as the finest topology on $X/N$ for which $\pi_{N}$
is continuous. Explicitly 
\begin{equation}
\tau_{N}:=\{A\subset X/N\mid\pi_{N}^{-1}(A)\in\tau\}.\label{qtop}
\end{equation}
Recall that $\tau_{N}$ is a linear topology on $X/N$ and that $\pi_{N}$
is linear, continuous, open, and onto. Moreover, $\tau_{N}$ is Hausdorff
separated if and only if $N$ is closed in $X$.

Let $(X,Y,\langle\cdot,\cdot\rangle)$ be a dual system. For $A\subset X$,
the \emph{orthogonal} \emph{of} $A$ is defined as: 
\begin{equation}
A^{\perp}:=\{y\in Y\mid\forall x\in A,\ \langle x,y\rangle=0\};\label{eq:-154}
\end{equation}
and $\sigma_{A}(y):=\iota_{A}^{*}(y)=\sup_{x\in A}\langle x,y\rangle$,
$y\in Y$ denotes the \emph{support functional }of $A$.

Similarly, for $B\subset Y$, $B^{\perp}:=\{x\in X\mid\forall y\in B,\ \langle x,y\rangle=0\}$
is the \emph{orthogonal} of $B$ and $\sigma_{B}(x)=\iota_{B}^{*}(x)=\sup_{y\in B}\langle x,y\rangle$,
$x\in X$ is the \emph{support functional of} $B$.

The pair $(\hat{X}:=X/Y^{\perp},\hat{Y}:=Y/X^{\perp})$ can be put
into duality via the well-defined bilinear form 
\begin{equation}
\langle\hat{x},\hat{y}\rangle:=\langle x,y\rangle,\ x\in X,\ y\in Y;\label{eq:-34-1}
\end{equation}
where $\hat{x}=x+Y^{\perp}$, $\hat{y}=y+X^{\perp}$. The dual system
$(X/Y^{\perp},Y/X^{\perp},\langle\cdot,\cdot\rangle)$ is separated
and it is called the \emph{separated version of} $(X,Y,\langle\cdot,\cdot\rangle)$.

For $T:X\rightrightarrows Y$, define $\hat{T}:X/Y^{\perp}\rightrightarrows Y/X^{\perp}$
the multivalued operator that makes the following diagram commutative,

\hspace{6cm}\begin{tikzcd}[row sep=3em, column sep=4em]
X 
  \arrow[d, "\pi_{Y^\perp}"'] 
  \arrow[r, shift left=0.6ex, "T"] 
  \arrow[r, shift right=0.6ex, swap, ""] 
  \arrow[dr, dashed] 
& Y 
  \arrow[d, "\pi_{X^\perp}"] \\
X/Y^{\perp} 
  \arrow[r, shift left=0.6ex, "\hat{T}"] 
  \arrow[r, shift right=0.6ex, swap, ""] 
& Y/X^{\perp}
\end{tikzcd}

\noindent i.e., $\pi_{X^{\perp}}\circ T=\hat{T}\circ\pi_{Y^{\perp}}$,
or equivalently 
\begin{equation}
(\hat{x},\hat{y})\in\operatorname*{Graph}(\hat{T})\Leftrightarrow\exists(u,v)\in\operatorname*{Graph}T,\ \hat{x}=\hat{u},\ \hat{y}=\hat{v}.\label{eq:-47}
\end{equation}
We call $\hat{T}$ the \emph{separated version of} $T$.

Similarly, for $A\subset X\times Y$, its \emph{separated version}
$\hat{A}$ is defined as 
\begin{equation}
(\hat{x},\hat{y})\in\hat{A}\Leftrightarrow\exists(u,v)\in A,\ \hat{x}=\hat{u},\ \hat{y}=\hat{v}.\label{eq:-17}
\end{equation}
Note that $\operatorname*{Graph}(\hat{T})$ is the separated version
of $\operatorname*{Graph}(T)$.

\begin{proposition} \label{ph} Let $(X,Y,\langle\cdot,\cdot\rangle)$
be a dual system, let $T:X\rightrightarrows Y$, and let $\hat{T}:\hat{X}\rightrightarrows\hat{Y}$
be the separated version of $T$. Then

\medskip{}

\emph{(i)} $T:X\rightrightarrows Y$ is monotone in $(X,Y,\langle\cdot,\cdot\rangle)$
if and only if $\hat{T}:\hat{X}\rightrightarrows\hat{Y}$ is monotone
in $(\hat{X},\hat{Y},\langle\cdot,\cdot\rangle)$.

\medskip{}

\emph{(ii)} $T:X\rightrightarrows Y$ is maximal monotone in $(X,Y,\langle\cdot,\cdot\rangle)$
if and only if $\hat{T}:\hat{X}\rightrightarrows\hat{Y}$ is maximal
monotone in $(\hat{X},\hat{Y},\langle\cdot,\cdot\rangle)$ and $\operatorname*{Graph}T$
is saturated with respect to $Y^{\perp}\times X^{\perp}$, that is,
\begin{equation}
\operatorname*{Graph}T=\operatorname*{Graph}T+Y^{\perp}\times X^{\perp}.\label{eq:}
\end{equation}
\end{proposition}
\begin{proof}
Subpoint (i) follows from $\langle x,y\rangle=\langle\hat{x},\hat{y}\rangle$,
$x\in X$, $y\in Y$. Note in addition that $(x,y)$ is monotonically
related to $T$ if and only if $(\hat{x},\hat{y})$ is monotonically
related to $\hat{T}$.

For (ii), assume first that $T$ is maximal monotone. If $(\hat{x},\hat{y})$
is monotonically related to $\hat{T}$, then any representative $(x,y)$
is monotonically related to $T$; hence, by maximality, $(x,y)\in\operatorname*{Graph}T$.
Thus $(\hat{x},\hat{y})\in\operatorname*{Graph}\hat{T}$, and $\hat{T}$
is maximal monotone. The same argument shows that $\operatorname*{Graph}T$
is saturated: if $(\hat{x},\hat{y})\in\operatorname*{Graph}\hat{T}$,
then every representative of $(\hat{x},\hat{y})$ is monotonically
related to $T$ and therefore belongs to $\operatorname*{Graph}T$.

Conversely, assume that $\hat{T}$ is maximal monotone and that $\operatorname*{Graph}T$
is saturated. If $(x,y)$ is monotonically related to $T$, then $(\hat{x},\hat{y})$
is monotonically related to $\hat{T}$, so $(\hat{x},\hat{y})\in\operatorname*{Graph}\hat{T}$,
that is, $(x,y)\in\operatorname*{Graph}T+Y^{\perp}\times X^{\perp}$.
By saturation, $(x,y)\in\operatorname*{Graph}T$. Therefore $T$ is
maximal monotone. 
\end{proof}
The natural duality $(Z:=X\times Y,Z,\cdot)$ can be similarly modified
to obtain its separated version $(\hat{Z}:=Z/Z^{\perp}=X/Y^{\perp}\times Y/X^{\perp},\hat{Z},\cdot)$,
where $Z^{\perp}=Y^{\perp}\times X^{\perp}$ and 
\begin{equation}
\hat{z}\cdot\hat{w}=z\cdot w,\ z,w\in Z.\label{eq:-27}
\end{equation}

Similarly, $A\subset X\times Y$ is monotone in $(X,Y,\langle\cdot,\cdot\rangle)$
if and only if $\hat{A}\subset\hat{Z}$ is monotone in $(\hat{Z},\hat{Z},\cdot)$.
Moreover, $A$ is maximal monotone if and only if $\hat{A}$ is maximal
monotone and $A$ is saturated with respect to $Z^{\perp}=Y^{\perp}\times X^{\perp}$,
i.e., $A=A+Z^{\perp}$. 

\begin{lemma} \label{P} Let $(X,Y,\langle\cdot,\cdot\rangle)$ be
a dual system and let $S,C\subset Z:=X\times Y$.

\medskip{}

\emph{(i)} $S$ is monotone (dissipative) with respect to the natural
duality $(Z,Z,\cdot)$ if and only if $\hat{S}$ is monotone (dissipative)
with respect to $(\hat{Z},\hat{Z},\cdot)$.

\medskip{}

\emph{(ii)} If $S$ is $\sigma(Z,C)-$compact then $\hat{S}$ is $\sigma(\hat{Z},\hat{C})-$compact.

\medskip{}

Here $\hat{S}:=\pi_{Z^{\perp}}(S)=\{\hat{s}=s+Z^{\perp}\mid s\in S\}$
is the separated version of $S$. \end{lemma}
\begin{proof}
(i) is due to $c(\hat{z})=c(z)$, $z\in Z$, and is a reformulation
of Proposition \ref{ph} (i); the dissipative case follows by applying
the monotone case to $\neg S$.

(ii) For every net $(z_{i})_{i}\subset Z$ such that $z_{i}\overset{\sigma(Z,C)}{\longrightarrow}z$
and every $\hat{c}=\pi_{Z^{\perp}}(c)\in\hat{C}$, with $c\in C$,
we have $\hat{z}_{i}\cdot\hat{c}=z_{i}\cdot c\overset{\mathbb{R}}{\to}z\cdot c$.
Therefore $\hat{z}_{i}\overset{\sigma(\hat{Z},\hat{C})}{\longrightarrow}\hat{z}$,
and so $\pi_{Z^{\perp}}:(Z,\sigma(Z,C))\to(\hat{Z},\sigma(\hat{Z},\hat{C}))$
is continuous. Hence, if $S$ is $\sigma(Z,C)-$compact, then $\hat{S}=\pi_{Z^{\perp}}(S)$
is $\sigma(\hat{Z},\hat{C})-$compact. 
\end{proof}
\begin{theorem} \label{Sep} Let $(X,Y,\langle\cdot,\cdot\rangle)$
be a dual system, let $G,S\subset X\times Y\eqqcolon Z$. Then 
\begin{equation}
G^{+}\cap S\neq\emptyset\ \Leftrightarrow\ \hat{G}^{+}\cap\hat{S}\neq\emptyset.\label{eq:-18}
\end{equation}
\end{theorem}
\begin{proof}
Assume that $z\in G^{+}\cap S$. Then $\hat{z}:=z+Z^{\perp}\in\hat{G}^{+}\cap\hat{S}\neq\emptyset$.

Conversely, let $\hat{z}\in\hat{G}^{+}\cap\hat{S}$. There is $s\in S$
such that $\hat{z}=\hat{s}$. For every $w\in G$, $\hat{w}\in\hat{G}$,
and since the quotient preserves the coupling, $c(s-w)=c(\hat{s}-\hat{w})=c(\hat{z}-\hat{w})\ge0$,
i.e., $s\in G^{+}\cap S\neq\emptyset$. 
\end{proof}
Similarly, for $T:X\rightrightarrows Y$, $S\subset X\times Y\eqqcolon Z$,
$\hat{S}=\pi_{Z^{\perp}}(S)$, and $\hat{T}:\hat{X}\rightrightarrows\hat{Y}$
defined as in (\ref{eq:-47}), 
\begin{equation}
\operatorname*{Graph}(T^{+})\cap S\neq\emptyset\ \Leftrightarrow\ \operatorname*{Graph}(\hat{T}^{+})\cap\hat{S}\neq\emptyset.\label{eq:-24}
\end{equation}

\section{General conditions for monotone extensions in dual systems}

Let $(X,Y,c=\langle\cdot,\cdot\rangle)$ be a dual system and let
$G,S\subset X\times Y$ be non-empty. Since $G^{+}=[\varphi_{G}\le c]$,
\begin{equation}
G^{+}\cap S=[\varphi_{G}-c+\iota_{S}\le0].\label{eq:-25}
\end{equation}

\noindent\tcbox[on line,colback=gray!15,colframe=black,arc=2pt,boxrule=0.4pt,left=3pt,right=3pt,top=2pt,bottom=2pt]{\textbf{Heuristic motivation.}}
Deciding whether $G^{+}\cap S$ is non-empty is the same as deciding
whether $\varphi_{G}-c$ takes a non-positive value on $S$. Equivalently,
it is the same as asking whether the function $f\coloneqq\varphi_{G}-c+\iota_{S}$
has a non-positive value. Therefore either the infimum of $f$ is
strictly negative, or the infimum of $f$ is zero and attained. More
precisely, $\inf_{Z}(\varphi_{G}-c+\iota_{S})=\inf_{S}(\varphi_{G}-c)$.
Hence 
\begin{equation}
G^{+}\cap S\neq\emptyset\ \Leftrightarrow\ \inf_{S}(\varphi_{G}-c)<0\ \mathrm{or}\ \min_{S}(\varphi_{G}-c)=0.\label{eq:-30}
\end{equation}
Here ``$\min$'' means that the infimum is finite and attained.

\medskip{}

We explain why dissipativity and convexity of $S$ are natural conditions
in this problem.

When $S$ is dissipative, the restriction $c|_{S}:S\to\mathbb{R}$,
$c|_{S}(z)\coloneqq c(z)$, is relatively $\sigma(Z,S)-$upper semicontinuous.
In other words, 
\begin{equation}
S\ \mathrm{dissipative}\ \Rightarrow\ \forall(s_{i})_{i}\subset S,\ s_{i}\overset{\sigma(Z,S)}{\to}s_{0}\in S:\ c(s_{0})\ge\limsup_{i}c(s_{i}).\label{eq:-32}
\end{equation}
Indeed, dissipativity gives $c(s_{i}-s_{0})=c(s_{i})+c(s_{0})-s_{i}\cdot s_{0}\le0$.
We find (\ref{eq:-32}) after we pass to the limit superior and take
into account that, since $s_{0}\in S$, $\lim_{i}s_{i}\cdot s_{0}=s_{0}\cdot s_{0}=2c(s_{0})$.

If $S$ is dissipative convex and $\sigma(Z,S)-$closed then $-c+\iota_{S}$
is convex $\sigma(Z,S)-$lower semicontinuous. Consequently, since
$\varphi_{G}$ is $\sigma(Z,G)-$lower semicontinuous, the function
$f=\varphi_{G}-c+\iota_{S}$ whose nonpositive values we need to detect
is itself convex and $\sigma(Z,G\cup S)-$lower semicontinuous.

Compactness of $S$ is natural for a complementary reason. The condition
$G^{+}\cap S\neq\emptyset$ is equivalent to asking whether $\varphi_{G}-c$
attains a non-positive value on $S$. Thus, after the preceding convexity
and semicontinuity considerations, compactness is the hypothesis that
permits limiting and minimization arguments on $S$. In particular,
when the infimum of $\varphi_{G}-c$ is equal to $0$, compactness
is what may turn this infimum into an attained minimum, and hence
into an actual point of $G^{+}\cap S$.

Let $z_{1},z_{2}\in S$. Assume that $z_{1},z_{2}\not\in G^{+}$,
that is, $(\varphi_{G}-c)(z_{1})>0$, $(\varphi_{G}-c)(z_{2})>0$.

For every $t\in[0,1]$, the convexity of $\varphi_{G}$ together with
the quadratic identity for $c$, gives the estimate 
\begin{equation}
(\varphi_{G}-c)(tz_{1}+(1-t)z_{2})\le q(t)\coloneqq t(\varphi_{G}-c)(z_{1})+(1-t)(\varphi_{G}-c)(z_{2})+t(1-t)c(z_{1}-z_{2}).\label{eq:-37}
\end{equation}
The quadratic upper bound $q(t)$ can take a non-positive value on
$[0,1]$ only if $c(z_{1}-z_{2})\le0$. Thus, for this mechanism to
work uniformly for pairs of points in $S$, it is natural to require
$S$ to be dissipative. Convexity of $S$ is also natural, because
once $z_{1},z_{2}\in S$, we want $tz_{1}+(1-t)z_{2}\in S$.

Under the preceding assumption, the quadratic upper bound $q(t)$
in (\ref{eq:-37}) has its minimizer inside $[0,1]$ precisely when
\[
c(z_{1}-z_{2})\le-|(\varphi_{G}-c)(z_{1})-(\varphi_{G}-c)(z_{2})|.
\]
The same quadratic function achieves a non-positive minimum on $[0,1]$
if and only if 
\[
c(z_{1}-z_{2})\le-(\sqrt{(\varphi_{G}-c)(z_{1})}+\sqrt{(\varphi_{G}-c)(z_{2})})^{2}.
\]

Consequently, dissipativity of $S$ is only a first-order qualitative
condition. To guarantee that $G^{+}\cap S$ is non-empty, one needs
additional hypotheses linking the geometry of $S$ with the size of
$\varphi_{G}-c$ on $S$.

\subsection{The canonical weak compact case}

The next result is the basic convex-analytic criterion behind the
extension theory. It expresses the nonemptiness of $G^{+}\cap S$
through a dual infimum condition. The subsequent proposition rewrites
this condition in more geometric forms, which then lead to directly
applicable criteria for monotone extension.

\begin{theorem} \label{E} Let $(X,Y,c)$ be a dual system, let $G\subset X\times Y$
be non-empty, and let $S\subset X\times Y$ be non-empty dissipative
convex and $\sigma(Z,G\cup S)-$compact. Then 
\begin{equation}
\min_{S}(\varphi_{G}-c)=\sup_{g\in\operatorname*{conv}G}\min_{s\in S}\left(s\cdot g-\operatorname*{conv}c_{G}(g)-c(s)\right);\label{eq:-23}
\end{equation}
\begin{equation}
\begin{aligned}G^{+}\cap S\neq\emptyset & \Leftrightarrow\ \min_{S}(\varphi_{G}-c)\le0\\
 & \Leftrightarrow\ \forall g\in\operatorname*{conv}G,\ \exists s\in S,\ (\operatorname*{conv}c_{G}-c)(g)+c(g-s)=\operatorname*{conv}c_{G}(g)+c(s)-s\cdot g\ge0.
\end{aligned}
\label{eq:-29}
\end{equation}
\end{theorem}
\begin{proof}
We first reduce to a Hausdorff separated setting on the linear span
that is relevant to the problem. Let $W\coloneqq\operatorname*{lin}(G\cup S)$
and endow $W$ with the restricted coupling. On $W$, the topology
induced by $\sigma(Z,G\cup S)$ is $\sigma(W,W)$. Passing, if necessary,
to the separated version of this restricted dual system and using
Lemma \ref{P} and Theorem \ref{Sep}, we may assume, without loss
of generality, that $Z=W$ and that $\sigma(Z,G\cup S)=\sigma(Z,Z)$
is Hausdorff separated.

If $\operatorname*{conv}c_{G}$ is not proper, then it takes the value
$-\infty$ at some point of $\operatorname*{conv}G$. Consequently
$\varphi_{G}=(\operatorname*{conv}c_{G})^{\square}\equiv+\infty$,
and both sides of (\ref{eq:-23}) are $+\infty$. In this case $G^{+}\cap S=\emptyset$,
and the last condition in (\ref{eq:-29}) fails at a point where $\operatorname*{conv}c_{G}=-\infty$.
Hence (\ref{eq:-29}) holds. We may therefore assume that $\operatorname*{conv}c_{G}$
is proper; then $\operatorname*{conv}c_{G}(g)\in\mathbb{R}$ for every
$g\in\operatorname*{conv}G$.

We use the Minimax Theorem for $f:S\times\operatorname*{conv}G\to\mathbb{R}$,
$f(s,g)\coloneqq\operatorname*{conv}c_{G}(g)+c(s)-s\cdot g$ for which,
for every $g\in\operatorname*{conv}G$, $f(\cdot,g)$ is concave and
$\sigma(Z,G\cup S)-$upper semicontinuous, while, for every $s\in S$,
$f(s,\cdot)$ is convex. We get 
\[
\inf_{g\in\operatorname*{conv}G}\max_{s\in S}f(s,g)=\max_{s\in S}\inf_{g\in\operatorname*{conv}G}f(s,g);
\]
see for example, \cite[Theorem 2]{MR47317} or \cite[Theorem 2.10.2]{MR1921556}.

Since $\inf_{g\in\operatorname*{conv}G}f(s,g)=c(s)-\varphi_{G}(s)$
for every $s\in S$, multiplying the minimax identity by $-1$ gives
(\ref{eq:-23}). The equivalences in (\ref{eq:-29}) follow from (\ref{eq:-25}),
(\ref{eq:-30}), and (\ref{eq:-23}). 
\end{proof}
\begin{remark} The set $G$ in Theorem \ref{E} is not assumed to
be monotone. Thus the result applies not only to monotone operator
graphs, but also to arbitrary comparison sets in $X\times Y$. The
monotonicity requirement enters only through the target condition
$G^{+}\cap S\neq\emptyset$, while the equivalent convex-analytic
conditions are expressed in terms of the Fitzpatrick-type function
$\varphi_{G}$ and the convexification $\operatorname*{conv}c_{G}$.
\end{remark}

\begin{proposition} \label{GR} Let $(X,Y,c=\langle\cdot,\cdot\rangle)$
be a dual system and let $G,S\subset X\times Y$. Consider the conditions

\medskip{}

\emph{(i)} \emph{$\forall g\in\operatorname*{conv}G,\ \exists s\in S,\ (\operatorname*{conv}c_{G}-c)(g)+c(g-s)\ge0$, }

\medskip{}

\emph{(ii)} \emph{$\forall g\in\operatorname*{conv}G,\ \exists s\in S,\ c(g-s)\ge0$, }

\medskip{}

\emph{(iii)} $\operatorname*{conv}(\Pr_{X}G)\subset\Pr_{X}S$.

\medskip{}

Then \emph{(iii) $\Rightarrow$ (ii).} If, in addition, $G$ is monotone,
then \emph{(ii) $\Rightarrow$ (i)}. \end{proposition}
\begin{proof}
(iii) $\Rightarrow$ (ii) When $\operatorname*{conv}(\Pr_{X}G)\subset\Pr_{X}S$,
for every $g\in\operatorname*{conv}G$ we take $s\in S$ with the
same first component as $g$; whence $c(g-s)=0$.

If, in addition, $G$ is monotone then (ii) $\Rightarrow$ (i) because
$\operatorname*{conv}c_{G}\ge c$. 
\end{proof}
\begin{theorem} \label{Em} Let $(X,Y,c)$ be a dual system, let
$G\subset X\times Y$ be non-empty monotone, and let $S\subset X\times Y$
be non-empty, dissipative, convex, and $\sigma(Z,G\cup S)-$compact.
Assume that 
\begin{equation}
\forall g\in\operatorname*{conv}G,\ \exists s\in S,\ c(g-s)\ge0.\label{eq:-28}
\end{equation}

The preceding condition is satisfied, in particular, if 
\begin{equation}
\operatorname*{Pr}\!_{X}G\subset\operatorname*{Pr}\!_{X}S.\label{eq:-6}
\end{equation}

Then $G^{+}\cap S\neq\emptyset$. \end{theorem}
\begin{proof}
Since $S$ is convex, $\operatorname*{Pr}\!_{X}G\subset\operatorname*{Pr}\!_{X}S$
implies $\operatorname*{conv}\operatorname*{Pr}\!_{X}G\subset\operatorname*{Pr}\!_{X}S$.

By Proposition \ref{GR} and Theorem \ref{E} we infer that $G^{+}\cap S\neq\emptyset$. 
\end{proof}
\begin{theorem} \label{ALT} Let $(X,Y,c)$ be a dual system, let
$G\subset X\times Y$ be non-empty, and let $S\subset X\times Y$
be non-empty dissipative convex and $\sigma(Z,G\cup S)-$compact.
Then exactly one of the following alternatives holds:

\medskip{}

\emph{(a)} $G$ admits a monotone extension into $S$, namely, $G^{+}\cap S\neq\emptyset$.

\medskip{}

\emph{(b)} There exist $g\in\operatorname*{conv}G$ and $\epsilon_{0}>0$
such that, for every $s\in S$, $(\operatorname*{conv}c_{G}-c)(g)+c(g-s)\le-\epsilon_{0}$.

\medskip{}

If, in addition, $G$ is monotone, and alternative \emph{(b)} holds,
then $S$ has a strong dissipative separation from a point of $\operatorname*{conv}G$;
that is, there exist $g\in\operatorname*{conv}G$ and $\epsilon_{0}>0$
such that, for every $s\in S$, $c(g-s)\le-\epsilon_{0}$. \end{theorem}
\begin{proof}
By (\ref{eq:-23}), (\ref{eq:-29}), either $\min_{S}(\varphi_{G}-c)\le0$,
which is equivalent to $G^{+}\cap S\neq\emptyset$, or 
\[
\sup_{g\in\operatorname*{conv}G}\min_{s\in S}\left(s\cdot g-\operatorname*{conv}c_{G}(g)-c(s)\right)>0,
\]
which is equivalent to (b). When $G$ is monotone we use $\operatorname*{conv}c_{G}\ge c$. 
\end{proof}
\begin{theorem} \label{H} Let $(X,Y,c=\langle\cdot,\cdot\rangle)$
be a dual system, let $G\subset X\times Y$ be non-empty, let $\omega\in Y$,
and let $K\subset X$ be non-empty, convex, and $\sigma(X,\operatorname*{Pr}_{Y}G\cup\{\omega\})-$compact.

Then $G^{+}\cap(K\times\{\omega\})\neq\emptyset$, i.e., there is
$\alpha\in K$ such that $(\alpha,\omega)$ is monotonically related
to $G$ if and only if 
\begin{equation}
\forall(x,y)\in X\times Y,\ \operatorname*{conv}c_{G}(x,y)\ge\inf_{\kappa\in K}\langle\kappa,y-\omega\rangle+\langle x,\omega\rangle,\label{eq:-26}
\end{equation}
if and only if 
\begin{equation}
\forall g\in\operatorname*{conv}G,\ \exists\kappa\in K,\ (\operatorname*{conv}c_{G}-c)(g)+c(g-(\kappa,\omega))\ge0.\label{eq:-31}
\end{equation}

Assume, in addition, that $G$ is monotone.

\medskip{}

\noindent\emph{(i)} If 
\[
\forall g\in\operatorname*{conv}G,\ \exists\kappa\in K,\ c(g-(\kappa,\omega))\ge0,
\]
then $G^{+}\cap(K\times\{\omega\})\neq\emptyset$.

\medskip{}

\noindent\emph{(ii)} If $\operatorname*{Pr}\!_{X}G\subset K$ then
$G^{+}\cap(K\times\{\omega\})\neq\emptyset$. In this case 
\begin{equation}
\forall(x,y)\in X\times Y,\ \operatorname*{conv}c_{G}(x,y)\ge\inf_{\kappa\in K}\langle\kappa,y-\omega\rangle+\langle x,\omega\rangle.\label{eq:-1}
\end{equation}
\end{theorem}
\begin{proof}
(\ref{eq:-26}) $\Leftrightarrow$ (\ref{eq:-31}) because $K$ is
$\sigma(X,\operatorname*{Pr}_{Y}G\cup\{\omega\})-$compact, in which
case, for every $y\in\operatorname*{conv}(\operatorname*{Pr}\!_{Y}G)$,
$\inf_{\kappa\in K}\langle\kappa,y-\omega\rangle$ is attained.

Set $S:=K\times\{\omega\}$. Then $S$ is non-empty, convex, and dissipative.
Moreover, the topology induced on $S$ by $\sigma(Z,G\cup S)$ is
weaker than the image of $\sigma(X,\operatorname*{Pr}_{Y}G\cup\{\omega\})$
under the map $\kappa\mapsto(\kappa,\omega)$; hence $S$ is $\sigma(Z,G\cup S)-$compact.
We use Theorem \ref{E} for this $S$ to conclude.

Assume, in addition, that $G$ is monotone. For (i) we use Theorem
\ref{Em} or (\ref{eq:-31}) and the inequality $\operatorname*{conv}c_{G}\ge c$,
which holds due to the monotonicity of $G$. Subpoint (ii) follows
as a direct consequence of Theorem \ref{Em} and (\ref{eq:-26}). 
\end{proof}
\begin{corollary} \label{HT} Let $(X,Y,c=\langle\cdot,\cdot\rangle)$
be a dual system, let $T:X\rightrightarrows Y$ be non-void, let $\omega\in Y$,
and let $K\subset X$ be non-empty, convex, and $\sigma(X,R(T)\cup\{\omega\})-$compact.

Then $\operatorname*{Graph}(T^{+})\cap(K\times\{\omega\})\neq\emptyset$,
i.e., there is $\alpha\in K$ such that $(\alpha,\omega)$ is monotonically
related to $T$, or equivalently, $\omega\in T^{+}(K)$ if and only
if 
\begin{equation}
\forall(x,y)\in X\times Y,\ \operatorname*{conv}c_{T}(x,y)\ge\inf_{\kappa\in K}\langle\kappa,y-\omega\rangle+\langle x,\omega\rangle,\label{eq:-33}
\end{equation}
if and only if 
\begin{equation}
\forall g\in\operatorname*{conv}(\operatorname*{Graph}T),\ \exists\kappa\in K,\ (\operatorname*{conv}c_{T}-c)(g)+c(g-(\kappa,\omega))\ge0.\label{eq:-34}
\end{equation}

Assume, in addition, that $T$ is monotone.

\medskip{}

\noindent\emph{(i)} If 
\[
\forall g\in\operatorname*{conv}(\operatorname*{Graph}T),\ \exists\kappa\in K,\ c(g-(\kappa,\omega))\ge0,
\]
then $\operatorname*{Graph}(T^{+})\cap(K\times\{\omega\})\neq\emptyset$.

\medskip{}

\noindent\emph{(ii)} If $D(T)\subset K$ then $\operatorname*{Graph}(T^{+})\cap(K\times\{\omega\})\neq\emptyset$.
In this case 
\begin{equation}
\forall(x,y)\in X\times Y,\ \operatorname*{conv}c_{T}(x,y)\ge\inf_{\kappa\in K}\langle\kappa,y-\omega\rangle+\langle x,\omega\rangle.\label{eq:-35}
\end{equation}
\end{corollary}
\begin{proof}
We apply Theorem \ref{H} for $G=\operatorname*{Graph}T$. 
\end{proof}
The next theorem is a compact-domain range principle. It shows that
compactness of a convex set containing the domain of a monotone operator
forces the range of its monotone polar to be large. In the maximal
monotone case, this becomes a range inclusion for the operator itself.

\begin{theorem} \label{GMC} Let $(X,Y,c=\langle\cdot,\cdot\rangle)$
be a dual system, let $T:X\rightrightarrows Y$ be non-void monotone,
let $C\subset Y$, and let $K\subset X$ be non-empty, convex, and
$\sigma(X,R(T)\cup C)-$compact.

If $D(T)\subset K$ then $C\subset T^{+}(K)$, or equivalently, for
every $\omega\in C$, $\operatorname*{Graph}(T^{+})\cap(K\times\{\omega\})\neq\emptyset$.

If, in addition, $T$ is maximal monotone then $C\subset R(T)$. \end{theorem}
\begin{proof}
Apply Corollary \ref{HT} (ii) for $S=K\times\{\omega\}$, for all
$\omega\in C$. 
\end{proof}
\begin{theorem} \label{GM} Let $(X,Y,c=\langle\cdot,\cdot\rangle)$
be a dual system, let $T:X\rightrightarrows Y$ be non-void monotone,
and let $K\subset X$ be non-empty, convex, and $\sigma(X,Y)-$compact.

If $D(T)\subset K$ then $T^{+}(K)=R(T^{+})=Y$.

If, in addition, $T$ is maximal monotone then $R(T)=Y$. \end{theorem}
\begin{proof}
Apply Theorem \ref{GMC} for $C=Y$. 
\end{proof}
\begin{corollary} \label{finit} Let $(X,Y,\langle\cdot,\cdot\rangle)$
be a dual system and let $G\subset X\times Y$ be finite monotone,
i.e., $G=\{(x_{i},y_{i})\}_{i\in\overline{1,n}}$ and, for every $i,j\in\overline{1,n}$,
$\langle x_{j}-x_{i},y_{j}-y_{i}\rangle\ge0$.

For every $y\in Y$ there is $x\in\operatorname*{conv}\{x_{i}\}_{i\in\overline{1,n}}$
such that $(x,y)$ is monotonically related to $G$, that is, 
\begin{equation}
\forall i\in\overline{1,n},\ \langle x-x_{i},y-y_{i}\rangle\ge0.\label{eq:-270}
\end{equation}
\end{corollary}
\begin{proof}
We use Theorem \ref{H} for the convex and $\sigma(X,Y)-$compact
$K:=\operatorname*{conv}\{x_{i}\}_{i\in\overline{1,n}}$ and $\omega=y$
after we note that $\operatorname*{Pr}\!_{X}G\subset K$. 
\end{proof}

\subsection{Compact target sets with auxiliary topologies}

The preceding results use compactness of $S$ with respect to the
canonical weak topology $\sigma(Z,G\cup S)$. In many applications,
however, the natural compactness may come from another topology on
$S$. The next result keeps the geometric assumptions that $S$ is
dissipative and convex, but replaces $\sigma(Z,G\cup S)-$compactness
by compactness with respect to an auxiliary Hausdorff topology $\beta$,
provided the coupling has the appropriate $\beta-$upper semicontinuity
properties on $S$.

\begin{theorem} \label{IM} Let $(X,Y,c)$ be a dual system, let
$G\subset X\times Y$ be non-empty, and let $S\subset X\times Y$
be non-empty, dissipative, and convex. Let $\beta$ be a Hausdorff
topology on $S$ such that $(S,\beta)$ is compact. Assume that 
\begin{equation}
\forall g\in\operatorname*{conv}G,\ \exists s\in S,\ (\operatorname*{conv}c_{G}-c)(g)+c(g-s)\ge0.\label{eq:-264}
\end{equation}
Assume further that

\medskip{}

\emph{(USC1)} For every $g\in\operatorname*{conv}G$, $c(g-\cdot):(S,\beta)\to\mathbb{R}$
is upper semicontinuous, in the sense that, for every $s_{0}\in S$
and for every net $\{s_{i}\}_{i\in I}\subset S$ such that $s_{i}\overset{\beta}{\to}s_{0}$,
we have $c(g-s_{0})\ge\limsup_{i\in I}c(g-s_{i})$.

\medskip{}

This condition is satisfied, in particular, under either of the following
stronger hypotheses:

\medskip{}

\emph{(USC2)} $c:(S,\beta)\to\mathbb{R}$ is upper semicontinuous
and $\beta$ is stronger than $\sigma(Z,G)$ on $S$; or

\medskip{}

\emph{(USC3)} $\beta$ is stronger than $\sigma(Z,G\cup S)$ on $S$.

\medskip{}

Then $G^{+}\cap S\neq\emptyset$.

\medskip{}

Moreover, \emph{(USC3)}$\Rightarrow$\emph{(USC2)}$\Rightarrow$\emph{(USC1)}.
\end{theorem}
\begin{proof}
First we show that (USC3)$\Rightarrow$(USC2)$\Rightarrow$(USC1).
Indeed, assume (USC3). Since $\beta$ is stronger than $\sigma(Z,S)$
on $S$ and $c:S\to\mathbb{R}$ is $\sigma(Z,S)-$upper semicontinuous
(recall (\ref{eq:-32})), $c$ is also $\beta-$upper semicontinuous.
Thus (USC2) holds.

Now assume (USC2). Let $z\in\operatorname*{conv}G$ and let $w_{i}\overset{\beta}{\to}w_{0}\in S$.
Since $c:(S,\beta)\to\mathbb{R}$ is upper semicontinuous, $c(w_{0})\ge\limsup_{i}c(w_{i})$.
Since $\beta$ is stronger than $\sigma(Z,G)$ on $S$, $\lim_{i}w_{i}\cdot z=w_{0}\cdot z$.
Using the identity $c(z-w_{i})=c(w_{i})+c(z)-w_{i}\cdot z$ we obtain
(USC1).

It remains to prove the result under (USC1).

For every $z\in\operatorname*{conv}G=\operatorname*{dom}(\operatorname*{conv}c_{G})$,
$\operatorname*{conv}c_{G}(z)\in\mathbb{R}$ due to (\ref{eq:-264}).

We use the Minimax Theorem for $A=(S,\beta)$, $B=\operatorname*{conv}G$,
and 
\[
f:S\times\operatorname*{conv}G\to\mathbb{R},\ f(w,z):=(\operatorname*{conv}c_{G}-c)(z)+c(z-w)=\operatorname*{conv}c_{G}(z)+c(w)-z\cdot w.
\]
For every $w\in S$, $f(w,\cdot):\operatorname*{conv}G\to\mathbb{R}$
is convex and $f(\cdot,z):S\to\mathbb{R}$ concave. More precisely,
for every $w_{1},w_{2}\in S$, $0\le t\le1$, we have $c(w_{1}-w_{2})\le0$
so 
\[
f(tw_{1}+(1-t)w_{2},z)=c(tw_{1}+(1-t)w_{2})+\operatorname*{conv}c_{G}(z)-(tw_{1}+(1-t)w_{2})\cdot z
\]
\[
=tf(w_{1},z)+(1-t)f(w_{2},z)-t(1-t)c(w_{1}-w_{2})\ge tf(w_{1},z)+(1-t)f(w_{2},z).
\]

From (USC1), for every fixed $z\in\operatorname*{conv}G$, $f(\cdot,z)=c(z-\cdot)+(\operatorname*{conv}c_{G}-c)(z)$
is $\beta-$upper semicontinuous. By the Minimax Theorem and (\ref{eq:-264})
\begin{equation}
\max_{w\in S}\inf_{z\in\operatorname*{conv}G}f(w,z)=\inf_{z\in\operatorname*{conv}G}\max_{w\in S}f(w,z)\ge0.\label{eq:-15}
\end{equation}

Notice that $\inf_{z\in\operatorname*{conv}G}f(w,z)=c(w)-(\operatorname*{conv}c_{G})^{\square}(w)=(c-\varphi_{G})(w)$,
$w\in Z$. Therefore 
\[
\max_{w\in S}\inf_{z\in\operatorname*{conv}G}f(w,z)=\max_{w\in S}(c-\varphi_{G})(w)\ge0,
\]
from which there exists $\overline{w}\in S$ such that $(c-\varphi_{G})(\overline{w})=\max_{w\in S}(c-\varphi_{G})(w)\ge0$,
i.e., $\overline{w}\in G^{+}\cap S$. 
\end{proof}
\begin{theorem} \label{IMm} Let $(X,Y,c)$ be a dual system, let
$G\subset X\times Y$ be non-empty monotone, and let $S\subset X\times Y$
be non-empty, dissipative, and convex. Let $\beta$ be a Hausdorff
topology on $S$ such that $(S,\beta)$ is compact.

Assume \emph{(USC1)}. This assumption is satisfied, in particular,
if either \emph{(USC2)} or \emph{(USC3)} holds, since \emph{(USC3)}$\Rightarrow$\emph{(USC2)}$\Rightarrow$\emph{(USC1)}.

Assume further that 
\begin{equation}
\forall z\in\operatorname*{conv}G,\ \exists w\in S,\ c(z-w)\ge0.\label{eq:-268}
\end{equation}
The preceding condition is satisfied, in particular, if 
\begin{equation}
\operatorname*{Pr}\,\!_{X}(G)\subset\operatorname*{Pr}\,\!_{X}(S).\label{eq:-269}
\end{equation}

Then $G^{+}\cap S\neq\emptyset$. \end{theorem}
\begin{proof}
As previously seen, (USC3)$\Rightarrow$(USC2)$\Rightarrow$(USC1)
and, since $G$ is monotone, (\ref{eq:-269}) $\Rightarrow$ (\ref{eq:-268})
$\Rightarrow$ (\ref{eq:-264}). Under (USC1), we use Theorem \ref{IM}
to get $G^{+}\cap S\neq\emptyset$. 
\end{proof}
\begin{remark} If $\beta$ is stronger than $\sigma(Z,G\cup S)$
on $S$, then compactness of $(S,\beta)$ implies compactness of $S$
with respect to $\sigma(Z,G\cup S)$. In this case, Theorem \ref{IM}
is covered by the canonical weak compact case, Theorem \ref{E}. Thus
the genuinely auxiliary-topology content of Theorem \ref{IM} appears
when $\beta$ is not stronger than $\sigma(Z,G\cup S)$ on $S$, while
still satisfying the upper semicontinuity requirement (USC1), for
instance through (USC2). The condition (USC3) is nevertheless useful
as a classical sufficient condition, since it implies (USC2), and
hence (USC1). \end{remark}

\begin{remark} \label{CC} Theorem \ref{IM} and Theorem \ref{IMm}
provide an alternative route to Theorem \ref{H} by applying them
with $S:=K\times\{\omega\}$ and $\beta$ equal to the trace topology
on $S$ of the product topology $\sigma(X,Y)\times\sigma(Y,\operatorname*{Pr}\,\!_{X}(G))$.
\end{remark}

\begin{theorem} \label{HIM} Let $(X,Y,c=\langle\cdot,\cdot\rangle)$
be a dual system, let $G\subset X\times Y$ be non-empty, let $K\subset X$
be non-empty and convex, and let $\omega\in Y$. Let $\beta$ be a
Hausdorff topology on $K$ such that $(K,\beta)$ is compact.

Assume that 
\begin{equation}
\forall(x,y)\in\operatorname*{conv}G,\ \exists u\in K,\ (\operatorname*{conv}c_{G}-c)(x,y)+\langle x-u,y-\omega\rangle\geq0.\label{eq:-19}
\end{equation}

Assume further that, for every $(a,b)\in\operatorname*{conv}G$, the
map $K\ni x\to\langle a-x,b-\omega\rangle$ is $\ensuremath{\beta}-$upper
semicontinuous. In particular, the upper semicontinuity assumption
is satisfied if either $K\ni x\to\langle x,\omega\rangle$ is $\ensuremath{\beta}-$upper
semicontinuous and $\beta$ is stronger than $\sigma(X,\operatorname*{Pr}_{Y}G)$
on $K$, or if $\beta$ is stronger than $\sigma\bigl(X,\operatorname*{Pr}_{Y}G\cup\{\omega\}\bigr)$
on $K$.

Then $G^{+}\cap(K\times\{\omega\})\neq\emptyset$, that is, there
exists $\alpha\in K$ such that $(\alpha,\omega)$ is monotonically
related to $G$. \end{theorem}
\begin{proof}
Apply Theorem \ref{IM} for $S:=K\times\{\omega\}$. 
\end{proof}
\begin{theorem} \label{HIMm} Let $(X,Y,c=\langle\cdot,\cdot\rangle)$
be a dual system, let $G\subset X\times Y$ be non-empty and monotone,
let $K\subset X$ be non-empty and convex, and let $\omega\in Y$.
Let $\beta$ be a Hausdorff topology on $K$ such that $(K,\beta)$
is compact.

Assume that 
\begin{equation}
\forall(x,y)\in\operatorname*{conv}G,\ \exists u\in K,\ \langle x-u,y-\omega\rangle\geq0.\label{eq:-20}
\end{equation}

The preceding condition is satisfied, in particular, if 
\begin{equation}
\operatorname*{Pr}\!_{X}G\subset K.\label{eq:-22}
\end{equation}

Assume further that, for every $(a,b)\in\operatorname*{conv}G$, the
map $K\ni x\to\langle a-x,b-\omega\rangle$ is $\ensuremath{\beta}-$upper
semicontinuous. In particular, the upper semicontinuity assumption
is satisfied if either $K\ni x\to\langle x,\omega\rangle$ is $\ensuremath{\beta}-$upper
semicontinuous and $\beta$ is stronger than $\sigma(X,\operatorname*{Pr}_{Y}G)$
on $K$, or if $\beta$ is stronger than $\sigma\bigl(X,\operatorname*{Pr}_{Y}G\cup\{\omega\}\bigr)$
on $K$.

Then $G^{+}\cap(K\times\{\omega\})\neq\emptyset$, that is, there
exists $\alpha\in K$ such that $(\alpha,\omega)$ is monotonically
related to $G$. \end{theorem}
\begin{proof}
Apply Theorem \ref{IMm} for $S:=K\times\{\omega\}$. 
\end{proof}
The preceding extension principles naturally contain variational inequalities
as a special case. Indeed, by choosing the target set to encode the
admissible values of the field and the comparison set to encode the
test points, the condition $G^{+}\cap S\neq\emptyset$ becomes precisely
the existence of a point satisfying the corresponding variational
inequality.

\begin{corollary} \label{VIS} Let $(X,Y,c=\langle\cdot,\cdot\rangle)$
be a dual system, let $C\subset X$ be non-empty, and let $S\subset X\times Y$
be non-empty, dissipative, and convex. Let $\beta$ be a Hausdorff
topology on $S$ such that $(S,\beta)$ is compact.

Assume that, for every $x\in\operatorname*{conv}C$, the map $S\ni(u,v)\to\langle u-x,v\rangle$
is $\ensuremath{\beta}-$upper semicontinuous.

Assume further that 
\begin{equation}
\forall x\in\operatorname*{conv}C,\ \exists(u,v)\in S,\ \langle u-x,v\rangle\geq0.\label{eq:-36}
\end{equation}

Then there exists $(a,b)\in S$ such that 
\[
\mathrm{VI}(S,C):\ \forall x\in\operatorname*{conv}C,\ \langle a-x,b\rangle\geq0.
\]
\end{corollary}
\begin{proof}
Apply Theorem \ref{IM} to $G:=C\times\{0\}$. Then $\operatorname*{conv}G=\operatorname*{conv}C\times\{0\}$,
$\ensuremath{\operatorname*{conv}c_{G}=0}$. Hence the geometric condition
in Theorem \ref{IM} becomes precisely (\ref{eq:-36}).

The upper semicontinuity condition in Theorem \ref{IM} is equivalent
to the $\ensuremath{\beta}-$upper semicontinuity of $S\ni(u,v)\to\langle u-x,v\rangle$,
for $x\in\operatorname*{conv}C$. 

Therefore $G^{+}\cap S\neq\emptyset$. Thus there exists $(a,b)\in S$
such that $(a,b)\in G^{+}$, namely, for every $x\in C$, $\langle a-x,b\rangle\geq0$. 
\end{proof}
The preceding results provide general criteria for dissipative compact
target sets $S\subset X\times Y$. We now specialize this framework
to target sets that are graphs and are not necessarily dissipative.
The next section treats a Browder-type situation, where the target
is the graph of a single-valued map defined on a compact convex set
containing the domain of the monotone operator. The following section
treats a Debrunner-Flor-type situation, where the target is again
the graph of a single-valued map, but the compact-domain containment
assumption is replaced by a convex-analytic compatibility condition.
The resulting graph-extension principles are then applied to variational
inequalities.

\section{An extension of Browder's theorem }

\begin{theorem} \label{BRS} Let $(X,Y,c=\langle\cdot,\cdot\rangle)$
be a dual system and let $K\subset X$ be a nonempty convex set. Let
$f:K\to Y$, and let $T:X\rightrightarrows Y$ be non-void monotone
such that $D(T)\subset K$. Let $\beta$ be a linear Hausdorff separated
topology on $\operatorname*{lin}K$ such that $(K,\beta)$ is compact.

Assume that, for every $z\in\operatorname{Graph}T$, the function
\begin{equation}
f_{z}:K\to\mathbb{R},\qquad f_{z}(x)\coloneqq c(z-(x,f(x)))\label{eq:-39}
\end{equation}
is $\beta-$upper semicontinuous, or, equivalently, whenever $\{x_{i}\}_{i\in I}\subset K$
is a net such that $x_{i}\overset{\beta}{\to}x_{0}$, one has $c(z-(x_{0},f(x_{0})))\ge\limsup_{i\in I}c(z-(x_{i},f(x_{i})))$.

Then $\operatorname*{Graph}(T^{+})\cap\operatorname*{Graph}(f)\neq\emptyset$.
Equivalently, there exists $x\in K$ such that $(x,f(x))$ is monotonically
related to $T$, that is, $f(x)\in T^{+}(x)$. \end{theorem}

\begin{remark} Define $\mathcal{T}_{\beta}$ as the strongest topology
on $\operatorname*{Graph}(f)$ that makes the map $E_{f}:(K,\beta)\to(\operatorname*{Graph}(f),\mathcal{T}_{\beta})$,
$E_{f}(x)=(x,f(x))$, continuous, i.e., 
\begin{equation}
\mathcal{T}_{\beta}:=\{U\subset\operatorname*{Graph}(f)\mid E_{f}^{-1}(U)=\operatorname*{Pr}\!_{X}U\in\beta\}.\label{eq:-255}
\end{equation}

The upper semicontinuity assumption in Theorem \ref{BRS} is implied
by the following more structural hypotheses. Suppose that $\mathcal{T}_{\beta}$
is stronger than $\sigma(Z,\operatorname*{Graph}T)$ on $\operatorname*{Graph}(f)$,
and that the coupling $c:(\operatorname*{Graph}(f),\mathcal{T}_{\beta})\to\mathbb{R}$
is upper semicontinuous. Equivalently, for every net $\{x_{i}\}_{i}\subset K$
with $x_{i}\overset{\beta}{\to}x_{0}$, and for every $(x,y)\in\operatorname*{Graph}T$
\begin{equation}
\lim_{i}\langle x_{i}-x_{0},y\rangle+\langle x,f(x_{i})-f(x_{0})\rangle=0,\ \limsup_{i}\langle x_{i},f(x_{i})\rangle\le\langle x_{0},f(x_{0})\rangle.\label{eq:-271}
\end{equation}
Then, for every $z\in\operatorname{Graph}T$, $f_{z}$ is $\beta-$upper
semicontinuous. \end{remark}

\begin{remark} Theorem \ref{BRS} generalizes Theorem \ref{BR},
which is obtained as a special case.

Indeed, assume the hypotheses of Theorem \ref{BR}. Apply Theorem
\ref{BRS} with $X=E$, $Y=F$, $T:X\rightrightarrows Y$ given by
$\operatorname*{Graph}T=G$, and $\beta=\tau$, which is a linear
Hausdorff separated topology on $\operatorname*{lin}K$ for which
$(K,\tau)$ is compact.

Since $f:(K,\tau)\to(F,\mu)$ is continuous, $E_{f}:(K,\tau)\to(\operatorname*{Graph}(f),\tau\times\mu)$,
$E_{f}(x)=(x,f(x))$ is continuous. Hence the topology $\mathcal{T}_{\tau}$
induced by $E_{f}$ on $\operatorname*{Graph}(f)$ is stronger than
the subspace topology inherited from $\tau\times\mu$.

By Remark \ref{RI}, equivalently from the continuity of the coupling
$c:(K\times F,\tau\times\mu)\to\mathbb{R}$, we know that, on $K\times F$,
$\tau\times\mu$ is stronger than $\sigma(Z=E\times F,\operatorname*{Graph}T)$.
Hence, on $\operatorname*{Graph}(f)\subset K\times F$, $\mathcal{T}_{\tau}$
is stronger than $\sigma(Z,\operatorname*{Graph}T)$. Moreover, the
continuity of $c:(K\times F,\tau\times\mu)\to\mathbb{R}$ implies
that $c:(\operatorname*{Graph}(f),\mathcal{T}_{\tau})\to\mathbb{R}$
is continuous. By the previous remark, for every $z\in\operatorname{Graph}T$,
$f_{z}$ is $\beta-$upper semicontinuous.

Thus all hypotheses in Theorem \ref{BRS} are satisfied. Applying
it with $\operatorname*{Graph}T=G$ gives exactly the conclusion of
Theorem \ref{BR}. Therefore Theorem \ref{BR} follows as a special
case of Theorem \ref{BRS}.

The converse implication is not available in general. The point is
that the new hypothesis of Theorem \ref{BRS} requires only the one-sided
upper semicontinuity of the functions $f_{z}$. This is substantially
weaker than requiring $f$ to be continuous into a weak topology on
$Y$, together with joint continuity of the coupling on a product
space. \end{remark}
\begin{proof}[Proof of Theorem \ref{BRS}]
Assume, by contradiction, that $\operatorname*{Graph}(T^{+})\cap\operatorname*{Graph}(f)=\emptyset$.
Hence, no element $w\in\operatorname*{Graph}(f)$ is monotonically
related to $T$, that is, 
\begin{equation}
\forall w\in\operatorname*{Graph}(f),\ \exists z\in\mathfrak{S}:=\operatorname*{Graph}T,\ c(z-w)<0,\label{eq:-256}
\end{equation}
or equivalently, since $\operatorname*{dom}f=K$, 
\begin{equation}
\forall x\in K,\ \exists z\in\mathfrak{S},\ c(z-(x,f(x)))<0.\label{eq:BRS-cover}
\end{equation}

For every $z\in\mathfrak{S}$, $f_{z}$ is $\beta-$upper semicontinuous.
Therefore 
\[
\forall z\in\mathfrak{S},\ O_{z}:=\{x\in K\mid f_{z}(x)=c(z-(x,f(x)))<0\}\ \mathrm{is}\ \beta-\mathrm{open\ in}\ K.
\]

Relation (\ref{eq:BRS-cover}) provides $K=\bigcup_{z\in\mathfrak{S}}O_{z}$.
Therefore, because $(K,\beta)$ is compact and Hausdorff, there exist
$m\ge1$, $\{z_{j}=(x_{j},y_{j})\}_{j\in\overline{1,m}}\subset\mathfrak{S}$,
such that $K=\bigcup_{j=1}^{m}O_{z_{j}}$ and $\{\rho_{1},\rho_{2},\ldots,\rho_{m}\}$
a partition of unity corresponding to this covering, i.e., for every
$j\in\overline{1,m}$, $\rho_{j}:K\to\mathbb{R}$ is $\beta-$continuous,
$\rho_{j}=0$ outside $O_{z_{j}}$, $0\le\rho_{j}\le1$, and $\sum_{j=1}^{m}\rho_{j}=1$. 

Consider $h:K\to X\times Y$, $h(x)=\sum_{j=1}^{m}\rho_{j}(x)z_{j}=(p(x),q(x))$,
where 
\[
p(x):=\sum_{j=1}^{m}\rho_{j}(x)x_{j},\ q(x):=\sum_{j=1}^{m}\rho_{j}(x)y_{j},\ x\in K.
\]
Because $\operatorname*{conv}D(T)\subset K$, for $K_{0}:=\operatorname*{conv}\{x_{1},\ldots,x_{m}\}$
the restriction $p|_{K_{0}}:K_{0}\to K_{0}$ is well-defined and $\beta-$con\-tin\-u\-ous.

The trace topology of $\beta$ on the finite-dimensional space $L:=\operatorname*{lin}\{x_{1},\ldots,x_{m}\}\subset\operatorname*{lin}K$
is the Euclidean topology because $\beta$ is linear and Hausdorff
separated. Since $K_{0}$ is ($\beta-$)com\-pact convex in $L$,
by the Brouwer fixed-point theorem, there is $x_{0}\in K_{0}$ such
that $p(x_{0})=x_{0}$.

Define $\lambda:K\to\mathbb{R}$ by 
\[
\lambda(x):=c(h(x)-(x,f(x)))=c(p(x)-x,q(x)-f(x))=c(\sum_{j=1}^{m}\rho_{j}(x)[z_{j}-(x,f(x))]).
\]
Then 
\[
\lambda(x)=\sum_{j=1}^{m}\rho_{j}(x)c(z_{j}-(x,f(x)))-\sum_{1\le i<j\le m}\rho_{i}(x)\rho_{j}(x)c(z_{i}-z_{j}),\ x\in K.
\]

Since $T$ is monotone, for every $x\in K$, $\sum_{1\le i<j\le m}\rho_{i}(x)\rho_{j}(x)c(z_{i}-z_{j})\ge0$.
Consequently, $\lambda(x)\le\sum_{j=1}^{m}\rho_{j}(x)c(z_{j}-(x,f(x)))$.
Because $\rho_{j}=0$ outside $O_{z_{j}}$, from the definition of
$O_{z_{j}}$ we get 
\begin{equation}
\forall x\in K,\ \forall j\in\overline{1,m},\ \rho_{j}(x)c(z_{j}-(x,f(x)))\le0;\label{eq:-257}
\end{equation}
from which we infer that 
\begin{equation}
\forall x\in K,\ \lambda(x)\le\sum_{j=1}^{m}\rho_{j}(x)c(z_{j}-(x,f(x)))\le0.\label{eq:-258}
\end{equation}

But $\lambda(x_{0})=0$. According to (\ref{eq:-257}), (\ref{eq:-258})
\begin{equation}
\forall j\in\overline{1,m},\ \rho_{j}(x_{0})c(z_{j}-(x_{0},f(x_{0})))=0.\label{eq:-259}
\end{equation}

Let $\mathscr{J}_{0}:=\{j\in\overline{1,m}\mid x_{0}\in O_{z_{j}}\}$.
Then $\mathscr{J}_{0}\neq\emptyset$ because $x_{0}\in K=\bigcup_{j=1}^{m}O_{z_{j}}$.
For every $j\in\mathscr{J}_{0}$, $x_{0}\in O_{z_{j}}$ from which,
according to (\ref{eq:-259}), $\rho_{j}(x_{0})=0$. For every $j\not\in\mathscr{J}_{0}$,
$\rho_{j}(x_{0})=0$ because $x_{0}\not\in O_{z_{j}}$. Therefore
$\sum_{j=1}^{m}\rho_{j}(x_{0})=\sum_{j\in\mathscr{J}_{0}}\rho_{j}(x_{0})+\sum_{j\not\in\mathscr{J}_{0}}\rho_{j}(x_{0})=0$
in contradiction to $\sum_{j=1}^{m}\rho_{j}=1$. 
\end{proof}
\begin{remark} Suppose one tries to recover Theorem \ref{BRS}, or
at least its particular case in which $\mathcal{T}_{\beta}$ is stronger
than $\sigma(Z,\operatorname*{Graph}T)$ on $\operatorname*{Graph}(f)$
and the coupling $c:(\operatorname*{Graph}(f),\mathcal{T}_{\beta})\to\mathbb{R}$
is upper semicontinuous, from Browder's Theorem \ref{BR}. Under the
assumptions of Theorem \ref{BRS}, if Theorem \ref{BR} is to be applied
with $\tau=\beta$, then the natural choice for the first ambient
space is $E=\operatorname*{lin}K$.

Assume in addition that $\beta$ is stronger than $\sigma(X,R(T))$
on $K$. In this case the first part in (\ref{eq:-271}) is equivalent
to the continuity of $f:(K,\beta)\to(Y,\sigma(Y,D(T)))$. This continuity
assumption aligns with the type considered in Theorem \ref{BR}, where
the second Browder space is $F=Y$ and $\mu=\sigma(Y,D(T))$. However,
in Theorem \ref{BR}, establishing the $\tau\times\mu-$continuity
of the coupling on $K\times Y$ is typically difficult.

From the point of view of $\beta$, one may replace the spaces $X$,
$Y$ by $\tilde{X}:=\operatorname*{lin}K$ and $\tilde{Y}:=\operatorname*{lin}R(T)$,
respectively. These spaces carry the restricted coupling $c$. Since
$\beta$ is Hausdorff separated, the dual system $(\tilde{X},\tilde{Y},c)$
is also Hausdorff in $\tilde{Y}$. If necessary, one may quotient
$\tilde{Y}$ by the annihilator of $\tilde{X}$, so that $(\tilde{X},\tilde{Y},c)$
becomes a separated dual system. Note also that $\sigma(\tilde{X},R(T))=\sigma(\tilde{X},\tilde{Y})$.

In this case, since $K$ is $\beta-$compact, it is also $\sigma(\tilde{X},\tilde{Y})-$compact.
Therefore $c:(K\times\tilde{Y},\sigma(\tilde{X},\tilde{Y})\times\mu(\tilde{Y},\tilde{X}))\to\mathbb{R}$
is continuous, where $\mu(\tilde{Y},\tilde{X})$ denotes the Mackey
topology compatible with the duality $(\tilde{X},\tilde{Y},c)$. This
is the strongest locally convex topology $\mu$ on $\tilde{Y}$ such
that the continuous dual $(\tilde{Y},\mu)^{*}=(\tilde{Y},\sigma(\tilde{X},\tilde{Y}))^{*}$.

Following this approach, the continuity of the coupling required in
Theorem \ref{BR} holds only if $\mu(\tilde{Y},\tilde{X})$ is weaker
than $\sigma(\tilde{Y},D(T))$, which is itself weaker than $\sigma(\tilde{Y},\tilde{X})$
because $D(T)\subset K$. Hence one would have to have $\mu(\tilde{Y},\tilde{X})=\sigma(\tilde{Y},\tilde{X})$
which can occur only when $\tilde{Y}$ and therefore $R(T)$, is finite-dimensional.

Therefore, when $R(T)$ is finite dimensional, this special case of
Theorem \ref{BRS} can be recovered from Theorem \ref{BR}. Outside
such finite-dimensional or otherwise exceptional situations, Theorem
\ref{BRS} is not a consequence of Theorem \ref{BR}.

The advantage of Theorem \ref{BRS} is precisely that it bypasses
joint continuity of the coupling and asks only for the upper semicontinuity
of the scalar functions $f_{z}$ in (\ref{eq:-39}). \end{remark}

\begin{theorem} \label{VI0} Let $(X,Y,\langle\cdot,\cdot\rangle)$
be a dual system. Let $K\subset X$ be convex, let $\emptyset\neq C\subset K$,
and let $F:K\to Y$. Let $\beta$ be a linear Hausdorff separated
topology on $\operatorname*{lin}K$ such that $(K,\beta)$ is compact.

Assume that, for every $a\in C$, $f_{a}:K\to\mathbb{R}$, $f_{a}(x)\coloneqq\langle a-x,F(x)\rangle$,
is $\beta-$upper semicontinuous. Equivalently, whenever $\{x_{i}\}_{i}\subset K$
is a net such that $x_{i}\overset{\beta}{\to}x_{0}$, one has 
\begin{equation}
\forall a\in C,\ \limsup_{i}\langle a-x_{i},F(x_{i})\rangle\le\langle a-x_{0},F(x_{0})\rangle.\label{eq:-8}
\end{equation}

Then there exists $x\in K$ that solves the variational inequality
\begin{equation}
\mathrm{VI}(F,K,C):\ \ \forall u\in\operatorname*{conv}C,\ \langle u-x,F(x)\rangle\ge0.\label{eq:-12}
\end{equation}

In particular, the upper semicontinuity condition (\ref{eq:-8}) holds
and $\mathrm{VI}(F,K,C)$ admits a solution if $F:(K,\beta)\to(Y,\sigma(Y,C))$
is continuous and, for every net $\{x_{i}\}_{i}\subset K$ with $x_{i}\overset{\beta}{\to}x_{0}$,
we have $\langle x_{0},F(x_{0})\rangle\le\liminf_{i}\langle x_{i},F(x_{i})\rangle$.
\end{theorem}
\begin{proof}
We use Theorem \ref{BRS} for $T:D(T)=C\subset K\rightrightarrows Y$,
$T(x)=\{0\}$, and $f=-F$. 
\end{proof}
As a first consequence, we recover the classical compactness-based
existence principle for Stampacchia variational inequalities, in a
form adapted to weakly compact sets and weak-to-strong continuous
fields.

\begin{corollary}[Weak compactness Stampacchia variational inequality]
Let $(H,\langle\cdot,\cdot\rangle_{H})$ be a real Hilbert space and
let $K\subset H$ be non-empty, convex, closed, and bounded. Let $F:K\to H$
be weak-to-strong continuous, that is, whenever a net $(x_{i})_{i}\subset K$
converges weakly to $x_{0}\in K$, the net $(F(x_{i}))_{i}$ converges
strongly to $F(x_{0})$. Then there exists $x\in K$ such that 
\[
\forall u\in K,\ \langle u-x,F(x)\rangle_{H}\geq0.
\]
\end{corollary}
\begin{proof}
Apply Theorem \ref{VI0} to the dual system $(H,H,\langle\cdot,\cdot\rangle_{H})$,
with $C=K$ and with $\beta$ equal to the weak topology on $H$.
Since $K$ is weakly compact, $(K,\beta)$ is compact.

Fix $a\in K$, and let $(x_{i})_{i}\subset K$ be a net such that
$x_{i}\to x_{0}$ weakly in $H$. By the weak-to-strong continuity
of $F$, we have $F(x_{i})\to F(x_{0})$ strongly in $H$. Hence 
\[
\langle a-x_{i},F(x_{i})\rangle_{H}=\langle a,F(x_{i})\rangle_{H}-\langle x_{i},F(x_{i})\rangle_{H}\to\langle a-x_{0},F(x_{0})\rangle_{H},
\]
because $\langle a,F(x_{i})\rangle_{H}\to\langle a,F(x_{0})\rangle_{H}$
and 
\[
\langle x_{i},F(x_{i})\rangle_{H}=\langle x_{0},F(x_{0})\rangle_{H}+\langle x_{i}-x_{0},F(x_{0})\rangle_{H}+\langle x_{i},F(x_{i})-F(x_{0})\rangle_{H}\to\langle x_{0},F(x_{0})\rangle_{H},
\]
since $x_{i}\to x_{0}$ weakly in $H$, $(x_{i})_{i}$ is weakly convergent
and therefore norm bounded, while $F(x_{i})-F(x_{0})\to0$ strongly
in $H$.

In particular, $x\mapsto\langle a-x,F(x)\rangle_{H}$ is weakly upper
semicontinuous on $K$. The conclusion follows from Theorem \ref{VI0}. 
\end{proof}

\section{An extension of the Debrunner-Flor theorem }

\begin{theorem} \label{DFE} Let $(X,Y,c=\langle\cdot,\cdot\rangle)$
be a dual system, let $K\subset X$ be a non-empty convex set, let
$f:K\to Y$, and let $T:X\rightrightarrows Y$ be non-void. Let $\tau$
be a Hausdorff separated locally convex topology on $X$ such that
$K$ is $\tau-$compact.

Assume that 
\begin{equation}
\forall\omega\in R(f),\ \forall g\in\operatorname*{conv}(\operatorname*{Graph}T),\ \exists\kappa\in K,\ (\operatorname*{conv}c_{T}-c)(g)+c(g-(\kappa,\omega))\ge0.\label{eq:-40}
\end{equation}

Assume further that, for every $g=(u,v)\in\operatorname*{conv}(\operatorname*{Graph}T)$,
the map 
\begin{equation}
K\times K\ni(x,x')\to c(g-(x',f(x)))=\langle x'-u,f(x)-v\rangle\label{eq:-41}
\end{equation}
is $\tau\times\tau-$upper semicontinuous.

Then $\operatorname*{Graph}(T^{+})\cap\operatorname*{Graph}(f)\neq\emptyset$,
that is, there exists $x\in K$ such that $(x,f(x))$ is monotonically
related to $T$, or equivalently $f(x)\in T^{+}(x)$.

Moreover, if $K$ is $\sigma(X,R(T)\cup R(f))-$compact, then the
geometric assumption (\ref{eq:-40}) above is equivalent to the infimum
condition 
\begin{equation}
\forall\omega\in R(f),\ \forall(x,y)\in\operatorname*{conv}(\operatorname*{Graph}T),\ \operatorname*{conv}c_{T}(x,y)\ge\inf_{\kappa\in K}\langle\kappa,y-\omega\rangle+\langle x,\omega\rangle.\label{eq:-42}
\end{equation}

In particular, the upper semicontinuity assumption (\ref{eq:-41})
is satisfied if there exists a topology $\mu$ on $R(f)$ such that
the restriction $c:(K\times R(f),\tau\times\mu)\to\mathbb{R}$ is
upper semicontinuous, $f:(K,\tau)\to(R(f),\mu)$ is continuous, $\tau$
is stronger than $\sigma(X,R(T))$ on $K$, and $\mu$ is stronger
than $\sigma(Y,D(T))$ on $R(f)$. \end{theorem}

\begin{remark} If $K$ is $\ensuremath{\tau}-$compact Hausdorff
and $\operatorname*{Graph}(f)$ is $\tau\times\mu-$compact in $K\times R(f)$,
then $f:(K,\tau)\to(R(f),\mu)$ is continuous. Indeed, the projection
$p:(\operatorname*{Graph}(f),\tau\times\mu)\to(K,\tau)$, $p(x,f(x))=x$,
is a continuous bijection from a compact space onto a Hausdorff space,
hence a homeomorphism. Then $f=\operatorname*{Pr}_{Y}\circ\left(p^{-1}\right)$
is continuous. Thus, under the compact Hausdorff assumption on $K$,
continuity of $f$ is equivalent to compactness of its graph. \end{remark}

\begin{remark} Theorem \ref{DF} is a special case of Theorem \ref{DFE}.
Indeed, take $K:=A$, let $f=\varphi:A\to Y$, and let $T:X\rightrightarrows Y$
be defined by $\operatorname*{Graph}T=M$. Then $T$ is monotone and
$D(T)\subset K$ by the assumptions of Theorem \ref{DF}.

The compactness and continuity hypotheses in Theorem \ref{DFE} are
also satisfied: $\tau$ is the given Hausdorff locally convex topology
of $X$, $K=A$ is $\ensuremath{\tau}-$compact, and $f=\varphi:A\to Y$
is continuous. Moreover, the assumed continuity of the bilinear form
on $X\times Y$ implies the upper semicontinuity condition (\ref{eq:-41})
appearing in Theorem \ref{DFE}. Finally, as observed above, the geometric
condition (\ref{eq:-40}) in Theorem \ref{DFE} follows from $D(T)\subset K$
and the monotonicity of $T$. \end{remark}
\begin{proof}
For $g=(u,v)\in\operatorname*{Graph}T$ and $w\in Y$ let 
\[
Q_{g}(w):=\{x\in X\mid(x,w)\ \mathrm{is\ monotonically\ related\ to}\ g,\ \mathrm{i.e.,}\ \langle x-u,w-v\rangle\ge0\}.
\]
Note that $x\in\bigcap_{g\in\operatorname*{Graph}T}Q_{g}(w)$ is equivalent
to $(x,w)$ is monotonically related to $T$. Define $M:Y\rightrightarrows X$,
\[
M(w):=K\cap\bigcap_{g\in\operatorname*{Graph}T}Q_{g}(w),\ w\in Y,
\]
and 
\[
F:=K\rightrightarrows K,\ F(x)=M(f(x)).
\]

With these notations, the conclusion is equivalent to the existence
of a fixed point for $F$, that is, there exists $x\in K$ such that
$x\in M(f(x))$, i.e., $x\in K$ and $(x,f(x))$ is monotonically
related to $T$.

We use the Ky Fan--Glicksberg fixed-point theorem. Since $K$ is
nonempty, convex, compact, and Hausdorff in the topology $\tau$,
it remains to verify that for every $x\in K$, $F(x)$ is a non-empty
convex subset of $K$ and that $\operatorname*{Graph}F$ is $\tau\times\tau-$closed
(see e.g. \cite{MR46638}).

We use Theorem \ref{HIM} for a fixed $\omega=f(x)\in R(f)$, where
$x\in K$, with $\beta=\tau$, $G\coloneqq\operatorname*{Graph}T$.
The upper semicontinuity condition in Theorem \ref{HIM} is obtained
from (\ref{eq:-41}) by fixing the first variable. By Theorem \ref{HIM}
we obtain that, for every $x\in K$, $\operatorname*{Graph}(T^{+})\cap\left(K\times\{f(x)\}\right)\neq\emptyset$,
that is, there is $\alpha\in K$ such that $(\alpha,f(x))$ is monotonically
related to $T$, so $\alpha\in K\cap\bigcap_{g\in\operatorname*{Graph}T}Q_{g}(f(x))=M(f(x))\neq\emptyset$,
from which $F$ has non-empty values.

The values of $F$ are convex because, for fixed $w\in Y$, each set
$Q_{g}(w)$ is convex and $K$ is convex.

It remains to prove that $\operatorname*{Graph}F$ is $\tau\times\tau-$closed.
Let $\{(x_{i},x_{i}')\}_{i\in I}\subset\operatorname*{Graph}F$ be
a net such that $x_{i}\overset{\tau}{\to}x$, $x_{i}'\overset{\tau}{\to}x'$.
Since $K$ is $\tau-$closed, $x,x'\in K$. For every $i\in I$, $x_{i}'\in F(x_{i})=(M\circ f)(x_{i})$.
Therefore $x_{i},x_{i}'\in K$ and $(x_{i}',f(x_{i}))$ is monotonically
related to $T$, that is, for every $(u,v)\in\operatorname*{Graph}T$,
\begin{equation}
\langle x_{i}'-u,f(x_{i})-v\rangle\ge0.\label{eq:-43}
\end{equation}

By the assumed upper semicontinuity (\ref{eq:-41}), passing to the
limit gives 
\begin{equation}
\forall(u,v)\in\operatorname*{Graph}T,\ \langle x'-u,f(x)-v\rangle\ge0.\label{eq:-44}
\end{equation}
Hence $x'\in F(x)$.

By the Ky Fan-Glicksberg fixed-point theorem, $F$ has a fixed point
$x\in K$, that is, $(x,f(x))\in\operatorname*{Graph}(T^{+})\cap\operatorname*{Graph}(f)$.

Now assume that $K$ is $\sigma(X,R(T)\cup R(f))-$compact. Then (\ref{eq:-40})
is equivalent to (\ref{eq:-42}) because, for every $y\in\operatorname*{conv}R(T)$
and every $\omega\in R(f)$, $\inf_{\kappa\in K}\langle\kappa,y-\omega\rangle$
is attained.

Finally, suppose that there exists a topology $\mu$ on $R(f)$ such
that the restriction $c:(K\times R(f),\tau\times\mu)\to\mathbb{R}$
is upper semicontinuous, $f:(K,\tau)\to(R(f),\mu)$ is continuous,
$\tau$ is stronger than $\sigma(X,R(T))$ on $K$, and $\mu$ is
stronger than $\sigma(Y,D(T))$ on $R(f)$.

Let $g=(u,v)\in\operatorname*{conv}(\operatorname*{Graph}T)$. Then
$u\in\operatorname*{conv}D(T)$, $v\in\operatorname*{conv}R(T)$.
For $x,x'\in K$, 
\[
\langle x'-u,f(x)-v\rangle=\langle u,v\rangle-\langle u,f(x)\rangle-\langle x',v\rangle+\langle x',f(x)\rangle.
\]
The term $(x,x')\to\langle x',f(x)\rangle$ is $\tau\times\tau-$upper
semicontinuous by the $\tau\times\mu-$upper semicontinuity of $c$
on $K\times R(f)$ and the continuity of $f$. The terms $x\to\langle u,f(x)\rangle$
and $x'\to\langle x',v\rangle$ are continuous by the assumptions
on $\mu$ and $\tau$ and the continuity of $f$.

Hence (\ref{eq:-41}) holds. The proof is complete. 
\end{proof}
\begin{theorem} \label{VID} Let $(X,Y,c=\langle\cdot,\cdot\rangle)$
be a dual system, let $C\subset X$ be non-empty, let $K\subset X$
be non-empty and convex, and let $F:K\to Y$. Let $\tau$ be a Hausdorff
separated locally convex topology on $X$ such that $K$ is $\ensuremath{\tau}-$compact.

Assume that 
\begin{equation}
\forall\omega\in R(F),\ \forall x\in\operatorname*{conv}C,\ \exists\kappa\in K,\ \langle\kappa-x,\omega\rangle\le0.\label{eq:-5}
\end{equation}

Assume further that, for every $a\in\operatorname*{conv}C$, the map
\begin{equation}
K\times K\ni(x,x')\to\langle a-x',F(x)\rangle\label{eq:-9}
\end{equation}
is $\ensuremath{\tau\times\tau}-$upper semicontinuous.

Then there exists $x\in K$ such that 
\[
\mathrm{VI}(F,K,C):\ \forall u\in\operatorname*{conv}C,\ \langle u-x,F(x)\rangle\geq0.
\]

In particular, the upper semicontinuity assumption above is satisfied
if there exists a topology $\mu$ on $R(F)$ such that $c:K\times R(F)\to\mathbb{R}$
is $\ensuremath{\tau\times\mu}-$lower semicontinuous, $F:(K,\tau)\to(R(F),\mu)$
is continuous, and $\mu$ is stronger than $\sigma(Y,C)$ on $R(F)$.
\end{theorem}
\begin{proof}
Apply Theorem \ref{DFE} for $T:D(T)=C\rightrightarrows Y$, $T(x)=\{0\}$,
and $f=-F$. In this case $\operatorname*{conv}(\operatorname*{Graph}T)=\operatorname*{conv}C\times\{0\}$,
$\operatorname*{conv}c_{T}=\iota_{\operatorname*{conv}C\times\{0\}}$,
(\ref{eq:-40}) becomes (\ref{eq:-5}), and (\ref{eq:-41}) transforms
into (\ref{eq:-9}). 
\end{proof}
The next corollary records the Hilbert-space form of the preceding
variational-inequality principle. In this setting, closed bounded
convex sets are weakly compact, and weak-to-strong continuity provides
the required scalar upper semicontinuity. 

\begin{corollary}[Weak compactness variational inequality] Let
$(H,\langle\cdot,\cdot\rangle_{H})$ be a real Hilbert space. Let
$C\subset H$ be non-empty, let $K\subset H$ be non-empty, convex,
closed and bounded, and let $F:K\to H$ be weak-to-strong continuous.
Assume that 
\begin{equation}
\forall\omega\in F(K),\ \forall x\in\operatorname*{conv}C,\ \exists\kappa\in K,\ \langle\kappa-x,\omega\rangle_{H}\leq0.\label{eq:-10}
\end{equation}

Then there exists $x\in K$ such that 
\[
\forall u\in\operatorname*{conv}C,\ \langle u-x,F(x)\rangle_{H}\geq0.
\]

In particular, when $C=K$, the geometric condition (\ref{eq:-10})
is automatically satisfied. Hence there exists $x\in K$ such that,
for every $u\in K$, $\langle u-x,F(x)\rangle_{H}\geq0$. \end{corollary}
\begin{proof}
Apply Theorem \ref{VID} with $X=Y=H$, with the Hilbert-space pairing,
and with $\tau$ equal to the weak topology on $H$. Since $K$ is
weakly compact, it is $\ensuremath{\tau}-$compact.

It remains to verify the upper semicontinuity condition. Let $a\in\operatorname*{conv}C$,
and suppose that $x_{i}\to x$ weakly in $K$ and $x_{i}'\to x'$
weakly in $K$. Since $F$ is weak-to-strong continuous, $F(x_{i})\to F(x)$
strongly. Hence $\langle a-x_{i}',F(x_{i})\rangle_{H}\to\langle a-x',F(x)\rangle_{H}$.

Thus the map $K\times K\ni(x,x')\to\langle a-x',F(x)\rangle_{H}$
is weak$\times$weak continuous, hence upper semicontinuous. Therefore
Theorem \ref{VID} applies. 
\end{proof}

\section{An auxiliary-topology example with a sharp threshold}

The purpose of the next example is to distinguish the hypotheses
of the extension principles. It satisfies the assumptions of
Theorem~\ref{IM}, but not the compactness assumption of
Theorem~\ref{E}. Its existence conclusion can also be obtained by
direct computation. 

\begin{example}
Let
\[
I=\{-1,0\}\cup\mathbb N,\qquad
X=\ell^{1}(I),\qquad Y=\ell^{\infty}(I),
\]
with the usual pairing $c(x,y)=\langle x,y\rangle=\sum_{i\in I}x_i y_i$, $x\in X, y\in Y$.

Write $\mathds{1}=(1)_{i\in I}$. Define

\begin{equation}
K=\left\{x\in X\mid 
x_i\ge0\ \text{for all }i\in I,\quad
x_{-1}=x_0=0,\quad \sum_{n\ge1}x_n\le1
\right\},
\end{equation}

and $m:X\to \mathbb R$, $\gamma:X\to X\times Y$

\begin{equation}
m(x)\coloneqq 
\sum_{n\ge1}x_n,\qquad
\gamma(x)\coloneqq \bigl(x,-m(x)\mathds1\bigr),\qquad 
S=\gamma (K)\subset X\times Y.
\end{equation}
Note that for $x\in K$ we have $m(x)=\langle x,\mathds1\rangle$. 

\medskip

\noindent {\bf Fact.1}
	The set $S$ is nonempty, convex, and dissipative.

\begin{proof}

	Indeed, $K$ is nonempty because $0\in K$, and it is convex because its
	defining conditions are preserved by convex combinations.
	Since $m$ and $\gamma$ are linear,  $S$ is nonempty and convex.
	
	The first two coordinates
	of $x,z\in K$ vanish. Hence  $\langle x-z,\mathds1\rangle=m(x)-m(z)$, $x,z\in K$.
	Consequently, $S$ is dissipative because, for every $x,z\in K$,

\[
c\bigl(\gamma (x)-\gamma (z)\bigr)=\left\langle
		x-z,-\bigl(m(x)-m(z)\bigr)\mathds1
		\right\rangle=-\bigl(m(x)-m(z)\bigr)^2\le0.
\]

\end{proof}	

Let $\tau=\sigma\bigl(\ell^{1}(I),c_{00}(I)\bigr)$ be the topology of coordinatewise convergence on
$X=\ell^{1}(I)$, where $c_{00}(I)\subset Y$ consists of the finitely supported
sequences. Since $\gamma |_K:K\to S$ is a bijection, define
a topology $\beta$ on $S$ by

\begin{equation}
U\in\beta
\quad\Longleftrightarrow\quad
(\gamma |_K)^{-1}(U)\text{ is relatively }\tau\text{-open in }K.	
\end{equation}

Thus $\gamma |_K:(K,\tau)\to(S,\beta)$ is a homeomorphism. For $(x^{(j)})_j\subset K$ and $x\in K$
\begin{equation}\label{cv}
\gamma (x^{(j)})\xrightarrow{\beta}\gamma (x)
\quad\Longleftrightarrow\quad
x^{(j)}_n\longrightarrow x_n
\quad\text{for every }n\ge1.
\end{equation}

\noindent {\bf Fact.2}
	The spaces $(K,\tau)$ and $(S,\beta)$ are compact and metrizable.

\begin{proof}
	Consider  $P=\{0\}^{\{-1,0\}}\times[0,1]^{\mathbb N}$ equipped with the product topology. Then $P$ is compact metrizable as a countable product of compact metrizable spaces. 
	
Inside $P$, we have
	\[
	K=\bigcap_{N\ge1}
	\left\{x\in P\mid\sum_{n=1}^{N}x_n\le1\right\}.
	\]
	Indeed, a nonnegative sequence whose partial sums are bounded
	by $1$ belongs to $\ell^{1}(I)$ and has total sum at most $1$.
	Each set in this intersection is closed in $P$, because finite
	coordinate sums are continuous. Hence $K$ is a closed subset
	of $P$, and therefore compact and metrizable.
	
The same properties hold for $(S,\beta)$, since $\gamma |_K$ is a homeomorphism.
\end{proof}

For every $i\in I$, let $e_i=(\delta_{ij})_{j\in I}$ denote
the standard coordinate vector, that is,
$(e_i)_j=0$, for $j\neq i$;  $(e_i)_i=1$.  We regard $e_i$ as an element of either $\ell^1(I)$ or
$\ell^\infty(I)$, according to context.

	For $\alpha>0$ and $n\ge1$, put
	
	\begin{equation}
u_n=\frac{\alpha}{2^{n+1}}(e_n-e_0)+n(e_{-1}-e_0), \qquad v_n=-2^{n+1}e_n,
\end{equation}
and define $ G_\alpha\coloneqq \{g_n=(u_n,v_n)\mid n\ge1\}\subset X\times Y$.
\end{example}
	
\noindent {\bf Fact.3}  $G_\alpha$ is dissipative and nonmonotone and both its projections are norm unbounded. 
	
\begin{proof}
Since $v_k$ has only one nonzero coordinate, at $k\ge1$,  $\langle u_n,v_k\rangle=-2^{k+1}(u_n)_k
=-\alpha\delta_{nk}$. 

Hence, for $n\neq k$, $\langle u_n-u_k,v_n-v_k\rangle=\langle u_n,v_n\rangle-\langle u_n,v_k\rangle
-\langle u_k,v_n\rangle+\langle u_k,v_k\rangle=-2\alpha<0$.

The three nonzero coordinates of $u_n$ occur at the distinct
indices $n,-1,0$. Thus $\|u_n\|_1
=\alpha/{2^{n+1}}+n+
\left(n+\alpha/{2^{n+1}}\right)
=2n+\alpha/{2^{n}}\to+\infty$. Therefore $\Pr_XG_\alpha$ is norm unbounded in $X$. Similarly, $\|v_n\|_\infty=2^{n+1}\to+\infty$, 
so $\Pr_YG_\alpha$ is norm unbounded in $Y$.

\end{proof}	
	
\noindent {\bf Fact.4}  	For $0<\alpha\le1$, all hypotheses of Theorem~\ref{IM} hold. More precisely \textup{(USC2)} holds and \textup{(USC3)} fails. 

\begin{proof}
Note first that $m:(K,\tau)\to\mathbb R$ is lower semicontinuous. Indeed, for every $N\ge1$, the function $K\ni x\longmapsto\sum_{n=1}^{N}x_n$ 
is $\tau$-continuous as a finite coordinate sum. 
Consequently, $m(x)
=\sup_{N\ge1}\sum_{n=1}^{N}x_n$ is $\tau$-lower semicontinuous. 
Since $m\ge0$ on $K$,
$-m^2$ is $\tau$-upper semicontinuous on $K$, because its strict sublevel sets are relatively $\tau $-open:  $[-m^2<a]=[m>\sqrt{-a}]$  , for $a\le0$; $[-m^2<a]=K$ for $a>0$.

Since $c(\gamma (x))=-m(x)^2$ and  $\gamma|_K:(K,\tau)\to(S,\beta)$ is
a homeomorphism, $c:(S,\beta)\to\mathbb R$ is upper semicontinuous, as required in
\textup{(USC2)}.

Moreover, $\beta=\sigma(Z,G_\alpha)|_S$.
Indeed, for $x\in K$ and $n\ge1$, using
$\langle u_n,\mathds1\rangle=0$, we obtain

\begin{equation}\label{tt}
\gamma(x)\cdot g_n=\langle x,v_n\rangle
	+\langle u_n,-m(x)\mathds1\rangle=-2^{n+1}x_n.
\end{equation}

Since each coordinate function $x\mapsto x_n$ is
$\tau$-continuous and
$\gamma|_K:(K,\tau)\to(S,\beta)$ is a homeomorphism,
every function $S\ni z\mapsto z\cdot g_n$ is
$\beta$-continuous. Therefore $\sigma(Z,G_\alpha)|_S\subseteq\beta$. The converse inclusion follows from (\ref{tt}), (\ref{cv}). 

Hence \textup{(USC2)} holds. However \textup{(USC3)} fails. Indeed,  $(e_n)_{n\ge1}\subset K$, $\gamma(e_1)\in S$, $\gamma(e_n)\xrightarrow{\beta}\gamma(0)=0$, while $\gamma(e_n)\cdot\gamma(e_1)=-2$, so $\gamma(e_n)$ does not converge to $0$ in $\sigma(Z,G_\alpha\cup  S)$. Thus $\beta$ is not stronger than
$\sigma(Z,G_\alpha\cup S)|_S$. 

Assume $0<\alpha\le1$. Fix $g=\sum_{n\in F}\lambda_ng_n\in\operatorname{conv}G_\alpha$, where, for all $n\in F$, $\lambda_n>0$ and $\sum_{n\in F}\lambda_n=1$. Here $F\subset\mathbb N$ is finite. 
Since for every $m\ge1$, 
$c(g_m)=-\alpha$, we know that $c_{G_\alpha}\ge-\alpha$. Hence  $\operatorname{conv}c_{G_\alpha}(g)=-\alpha$.

Using $c(g-z)=c(g)+c(z)-g\cdot z$, together with
\[
c(\gamma(x))=-m(x)^2,\qquad
\gamma(x)\cdot g_n=-2^{n+1}x_n,
\]
we obtain, for every $x\in K$,
\[
(\operatorname{conv}c_{G_\alpha}-c)(g)+c(g-\gamma(x))
=-\alpha-m(x)^2+\sum_{n\in F}\lambda_n2^{n+1}x_n.
\]

There exists $j\in F$ such that $\lambda_j2^{j+1}>2$.
Otherwise, for every $j\in F$, $\lambda_j\le2^{-j}$ from which 
\[
1=\sum_{n\in F}\lambda_n
\le\sum_{n\in F}2^{-n}
<\sum_{n=1}^{\infty}2^{-n}=1,
\]
a contradiction.
Taking $x=e_j\in K$, for which $m(e_j)=1$, gives
\[
(\operatorname{conv}c_{G_\alpha}-c)(g)
	+c(g-\gamma(e_j))=-\alpha-1+\lambda_j2^{j+1}
	>1-\alpha\ge0.
\]
Since $\gamma(e_j)\in S$, the geometric condition
of Theorem~\ref{IM} holds.

	Thus Theorem~\ref{IM} applies and yields
\[
G_\alpha^+\cap S\neq\varnothing
\qquad(0<\alpha\le1).
\]
\end{proof}

\noindent {\bf Fact.5} For every $\alpha>0$

\begin{equation}\label{min}
	\min_{s\in S}(\varphi_{G_\alpha}-c)(s)=\alpha-1,
\end{equation}
and therefore

\begin{equation}
	G_\alpha^+\cap S\neq\varnothing
	\quad\Longleftrightarrow\quad 0<\alpha\le1.
\end{equation}

At $\alpha=1$, this intersection is the singleton
\[
G_1^+\cap S=\{(\bar x,-\mathds1)\},
\]
where
\[
\bar x_{-1}=\bar x_0=0,\qquad
\bar x_n=2^{-n}\quad(n\ge1).	
\]

\begin{proof}
For every $x\in K$,
$\varphi_{G_\alpha}(\gamma(x))=\sup_{n\ge1}\bigl(\gamma(x)\cdot g_n-c(g_n)\bigr)=\alpha-\inf_{n\ge1}2^{n+1}x_n$.

Since $c(\gamma(x))=-m(x)^2$, we obtain
\[
(\varphi_{G_\alpha}-c)(\gamma(x))
=\alpha+m(x)^2-\inf_{n\ge1}2^{n+1}x_n.
\]
Set $h=\inf_{n\ge1}2^{n+1}x_n$. Then $h\ge0$ and
$x_n\ge2^{-(n+1)}h$ for every $n\ge1$, whence
\[
m(x)\ge h\sum_{n\ge1}2^{-(n+1)}=\frac h2.
\]
Consequently,

\begin{equation}\label{est}
(\varphi_{G_\alpha}-c)(\gamma(x))\ge\alpha+m(x)^2-2m(x)=\alpha-1+(m(x)-1)^2\ge\alpha-1.
\end{equation}
Define $\bar x\in K$ by
\[
\bar x_{-1}=\bar x_0=0,\qquad
\bar x_n=2^{-n}\quad(n\ge1).
\]
Then $m(\bar x)=1$ and $2^{n+1}\bar x_n=2$ for every
$n\ge1$, so the preceding lower bound is attained at
$\gamma(\bar x)$. Thus
\[
\min_{z\in S}(\varphi_{G_\alpha}-c)(z)=\alpha-1.
\]
Since $G_\alpha^+=[\varphi_{G_\alpha}\le c]$, it follows that,
for $\alpha>0$,
\[
G_\alpha^+\cap S\neq\varnothing
\quad\Longleftrightarrow\quad \alpha\le1.
\]

If $\alpha=1$ and $\gamma(x)\in G_1^+\cap S$, for some $x\in K$, the preceding
estimate (\ref{est}) implies $m(x)=1$. The inequality
$(\varphi_{G_1}-c)(\gamma(x))\le0$ then gives
\[
\inf_{n\ge1}2^{n+1}x_n\ge2,
\]
so $x_n\ge2^{-n}$ for every $n\ge1$. Since
$
\sum_{n\ge1}x_n=1=\sum_{n\ge1}2^{-n}$,
we conclude that $x=\bar x$. Therefore
\[
G_1^+\cap S
=\{\gamma(\bar x)\}
=\{(\bar x,-\mathds{1})\}.
\]
\end{proof}

\noindent {\bf Fact.6}  $S$ is not $\sigma(Z,G_\alpha\cup S)$-compact. In particular, this noncompactness shows again that (USC3) fails. 

\begin{proof}
Suppose that $S$ is $\sigma(Z,G_\alpha\cup S)$-compact.
Then the sequence $(\gamma(e_n))_{n\ge1}\subset S$ has
a subnet $(\gamma(e_{n_\iota}))_\iota$ converging in
$\sigma(Z,G_\alpha\cup S)$ to some $\gamma(x)\in S$.

For every fixed $j\ge1$, since $n_\iota\to+\infty$,
\[
\gamma(e_{n_\iota})\cdot g_j
=-2^{j+1}\delta_{n_\iota j}\longrightarrow0.
\]
Hence
\[
\gamma(x)\cdot g_j=-2^{j+1}x_j=0
\qquad(j\ge1).
\]
Since $x\in K$, this implies $x=0$.

On the other hand, $\gamma(e_1)\in S$ and
\[
\gamma(e_{n_\iota})\cdot\gamma(e_1)=-2
\qquad\text{for every }\iota,
\]
whereas $\gamma(0)\cdot\gamma(e_1)=0$.
This contradicts the assumed convergence.
Therefore $S$ is not $\sigma(Z,G_\alpha\cup S)$-compact.
\end{proof}

\begin{remark}
	This example satisfies \textup{(USC2)} but not
	\textup{(USC3)}, and lies outside the compactness assumptions
	of Theorem~\ref{E}. The classical extension theorems requiring
	monotonicity of the comparison set do not apply to $G_\alpha$.
	Here the conclusion concerns a point monotonically related
	to $G_\alpha$; adjoining that point does not make the already
	nonmonotone set $G_\alpha$ monotone.
\end{remark}

\section{A compact-target criterion for global minimization}

The next application explores the problem of finding a global minimizer of a potentially nonconvex function inside a prescribed compact set.

\begin{proposition}
	Let $f:\mathbb R^d\to\mathbb R$ be continuously differentiable,
	and let $K\subset\mathbb R^d$ be nonempty, compact, and convex.
	Set
	\[
	G=\operatorname{Graph}(\nabla f),\qquad
	c(u,v)=\langle u,v\rangle,\ (u,v)\in\mathbb R^d\times\mathbb R^d,\qquad
	c_G=c+\iota_G.
	\]
	Assume that
	\begin{equation}\label{opt-criterion}
		\operatorname{conv}c_G(u,v)
		\ge\inf_{\kappa\in K}\langle\kappa,v\rangle
		\qquad((u,v)\in\mathbb R^d\times\mathbb R^d).
	\end{equation}
	Then there exists $\bar x\in K$ such that
	\begin{equation}\label{opt-minty}
		\langle u-\bar x,\nabla f(u)\rangle\ge0
		\qquad(u\in\mathbb R^d).
	\end{equation}
	In particular, $f(\bar x)=\min_{u\in\mathbb R^d}f(u)$. 	Moreover, $f$ is nondecreasing along every ray issuing from
	$\bar x$.

\end{proposition}

\begin{proof}
	Apply Theorem~\ref{H} to the Euclidean dual system with
	\[
	G=\operatorname{Graph}(\nabla f),\qquad \omega=0.
	\]
	The compactness assumption holds because $K$ is compact
	in the Euclidean topology. Condition~\eqref{opt-criterion}
	is precisely the convex-analytic criterion of that theorem.
	Hence there exists $\bar x\in K$ such that $(\bar x,0)\in G^+\cap(K\times\{0\})$. Monotone relatedness
	gives
	\[
	\langle\bar x-u,-\nabla f(u)\rangle\ge0
	\qquad(u\in\mathbb R^d),
	\]
	which is \eqref{opt-minty}.
	
	Fix $y\in\mathbb R^d$ and define
	\[
	h(t)=f\bigl(\bar x+t(y-\bar x)\bigr),\qquad t\ge0.
	\]
	For every $t>0$, applying \eqref{opt-minty} at
	$u=\bar x+t(y-\bar x)$ yields
	\[
	t\,h'(t)
	=\left\langle t(y-\bar x),
	\nabla f\bigl(\bar x+t(y-\bar x)\bigr)\right\rangle
	\ge0.
	\]
	Thus $h$ is nondecreasing on $[0,+\infty)$, and
	\[
	f(y)=h(1)\ge h(0)=f(\bar x).
	\]
	Since $y$ is arbitrary, $\bar x$ is a global minimizer.
\end{proof}

\begin{remark}
	By Theorem~\ref{H}, condition~\eqref{opt-criterion} is
	necessary and sufficient for the existence of a point
	$\bar x\in K$ satisfying \eqref{opt-minty}.
	It is sufficient, but need not be necessary, for the
	existence of a global minimizer in $K$.
\end{remark}

\begin{proposition}[A diameter criterion for the Cauchy location objective]
	Let $b_1,\ldots,b_N\in\mathbb R^d$, let $w_1,\ldots,w_N>0$,
	and let $\sigma>0$. Set
	\[
	f(x)=\sum_{j=1}^N w_j\log\bigl(\sigma^2+\|x-b_j\|^2\bigr),\ x\in {\mathbb R}^d,
	\quad K=\operatorname{conv}\{b_1,\ldots,b_N\},
	\quad D=\max_{i,j}\|b_i-b_j\|.
	\]
	Let $G=\operatorname{Graph}(\nabla f)$. If $D\le 2\sigma$, then
	\[
	\operatorname{conv}c_G(u,v)\ge\min_{k\in K}\langle k,v\rangle
	\qquad\bigl((u,v)\in\mathbb R^d\times\mathbb R^d\bigr).
	\]
	Consequently, there exists $\bar x\in K$ such that
	\begin{equation}\label{eq:cauchy-minty}
		\langle u-\bar x,\nabla f(u)\rangle\ge0
		\qquad\text{for every }u\in\mathbb R^d.
	\end{equation}
	Moreover, $f$ is strictly increasing along every ray issuing from
	$\bar x$. In particular, $\bar x$ is its unique global minimizer.
	The constant $2$ is sharp as a uniform diameter bound.
\end{proposition}

\begin{proof}
	Write $r_j(x)=\sigma^2+\|x-b_j\|^2$ and set
	\[
	A(x)=\sum_{j=1}^N\frac{w_j}{r_j(x)},\qquad
	p_j(x)=\frac{w_j/r_j(x)}{A(x)},\qquad
	T(x)=\sum_{j=1}^N p_j(x)b_j.
	\]
	Since the $p_j(x)$ are positive and sum to one, $T(x)\in K$.
	Moreover,
	\[
	\nabla f(x)=2A(x)(x-T(x)).
	\]
	
	Fix $x$ and regard $j$ as a random index with probabilities $p_j(x)$.
	For unit vectors $h,\ell$, differentiation gives
	\[
	\langle\ell,DT(x)h\rangle
	=\operatorname{Cov}_{p(x)}\left(
	\langle\ell,b_j\rangle,
	\frac{2\langle b_j-x,h\rangle}{r_j(x)}\right).
	\]
	
	Here $DT(x)h=\left.\frac{d}{dt}T(x+th)\right|_{t=0}$ is the derivative of $T:\mathbb R^d\to\mathbb R^d$ at $x$
	applied to $h$, while  $\operatorname{Cov}_{p(x)}$ stands for the covariance with respect to the weights \(p_1(x),\ldots,p_N(x)\). 
	
	\eject

	The first variable lies in an interval of length at most $D$.
	Its variance is therefore at most $D^2/4$, as is seen by
	centering at the midpoint of that interval.
	The second variable has absolute value at most $1/\sigma$.
	By the Cauchy--Schwarz inequality,
	$\|DT(x)\|\le D/(2\sigma)\le1$.
	Hence $T$ is nonexpansive and $M=I-T$ is monotone.
	
	Take any finite convex representation
	\[
	(u,v)=\sum_i\lambda_i(x_i,\nabla f(x_i)).
	\]
	Put $a_i=2A(x_i)$, $a=\sum_i\lambda_i a_i$, and
	$\mu_i=\lambda_i a_i/a$. Set
	\[
	y=\sum_i\mu_i x_i,\qquad k=\sum_i\mu_iT(x_i)\in K.
	\]
	Since $\sum_i\mu_iM(x_i)=y-k$ and $\sum_i\mu_i(x_i-y)=0$,
	\[
	\begin{aligned}
		\frac1a\sum_i\lambda_i\langle x_i-k,\nabla f(x_i)\rangle
		&=\sum_i\mu_i\langle x_i-y,M(x_i)-M(y)\rangle
		+\|y-k\|^2\\
		&\ge0.
	\end{aligned}
	\]
	It follows that
	\[
	\sum_i\lambda_i\langle x_i,\nabla f(x_i)\rangle
	\ge\langle k,v\rangle\ge\min_{z\in K}\langle z,v\rangle.
	\]
	Taking the infimum over all such representations yields the
	claimed inequality on $\operatorname{conv}G$. Outside
	$\operatorname{conv}G$ its left-hand side is $+\infty$.
	
	Theorem~13 with $\omega=0$ now yields $\bar x\in K$ satisfying
	\eqref{eq:cauchy-minty}. For $h\ne0$, put $q(t)=f(\bar x+th)$.
	Then $tq'(t)\ge0$ for $t>0$, whence $q$ is nondecreasing on
	$[0,\infty)$. Equality at two distinct points would make $q$
	constant on the interval between them, hence everywhere by real
	analyticity. This contradicts $q(t)\to+\infty$ as $t\to+\infty$.
	Thus $q$ is strictly increasing on $[0,\infty)$.
	
	For sharpness, take two observations $-b,b$ of unit weight in
	$\mathbb R$, where $b>\sigma$. Then $f$ is even and
	\[
	f''(0)=\frac{4(\sigma^2-b^2)}{(\sigma^2+b^2)^2}<0.
	\]
	For all sufficiently small $t>0$, one has $f'(t)<0$.
	The graph points $(t,f'(t))$ and $(-t,-f'(t))$ have midpoint
	$(0,0)$. Therefore
	\[
	\operatorname{conv}c_G(0,0)\le t f'(t)<0
	=\min_{k\in K'}\langle k,0\rangle
	\]
	for every nonempty compact convex $K'\subset\mathbb R$.
\end{proof}

\begin{remark}
	The objective can be nonconvex even on $K$.
	Indeed, take $\sigma=1$, $b_1=0$, $b_2=\sqrt3$, $w_1=9$,
	and $w_2=1$. Then $D=\sqrt3<2$, whereas
	\[
	f(x)=9\log(1+x^2)+\log\bigl(1+(x-\sqrt3)^2\bigr)
	\quad\text{and}\quad f''(\sqrt3)=-\tfrac14.
	\]
	Thus $f''<0$ at interior points of $K=[0,\sqrt3]$ sufficiently
	close to $\sqrt3$.
\end{remark}

\begin{proposition}[A location estimate with a variable reference point]
	Let $K\subset\mathbb R^d$ be compact and convex, let
	$b_1,\ldots,b_N\in K$, and let $w_1,\ldots,w_N,\sigma,\lambda>0$.
	Assume that $\operatorname{diam}K\le 2\sigma$. For $p\in K$, put
	\[
	\Phi_p(u)=\sum_{j=1}^N w_j
	\log\bigl(\sigma^2+\|u-b_j\|^2\bigr)
	+\frac{\lambda}{2}\|u-p\|^2.
	\]
	Then, for every continuous $q:K\to K$, there exists $\bar x\in K$
	such that
	\[
	\operatorname*{argmin}_{u\in\mathbb R^d}\Phi_{q(\bar x)}(u)
	=\{\bar x\}.
	\]
	Moreover, $\Phi_{q(\bar x)}$ is strictly increasing along every
	ray issuing from $\bar x$.
\end{proposition}

\begin{proof}
	Set
	\[
	d_j(u)=\frac{2w_j}{\sigma^2+\|u-b_j\|^2},\qquad
	H(u)=\lambda u+\sum_j d_j(u)(u-b_j),\qquad f(x)=\lambda q(x).
	\]
	For $p\in K$, let
	\[
	a(u)=\lambda+\sum_jd_j(u),\qquad
	T_p(u)=\frac{\lambda p+\sum_jd_j(u)b_j}{a(u)}.
	\]
	Then $T_p(u)\in K$ and
	\[
	\nabla\Phi_p(u)=H(u)-\lambda p=a(u)(u-T_p(u)).
	\]
	The derivative estimate for the Cauchy location objective applies
	to $T_p$ with one additional observation at $p$ of constant weight
	$\lambda$. Its logarithmic derivative is zero, and all observations
	belong to $K$. Thus
	\[
	\|DT_p(u)\|\le\frac{\operatorname{diam}K}{2\sigma}\le1.
	\]
	It follows that $I-T_p$ is monotone.
	
	Fix $p\in K$. By coercivity, $\Phi_p$ has a global minimizer $y_p$.
	Since $\nabla\Phi_p(y_p)=0$, we have $y_p=T_p(y_p)\in K$. Hence
	\[
	\langle u-y_p,H(u)-\lambda p\rangle
	=a(u)\langle u-y_p,(I-T_p)(u)-(I-T_p)(y_p)\rangle\ge0
	\qquad(u\in\mathbb R^d).
	\]
	Consequently, the affine function
	\[
	(u,v)\longmapsto
	\langle u,\lambda p\rangle+\langle y_p,v-\lambda p\rangle
	\]
	minorizes $c$ on $\operatorname{gra}H$, and hence minorizes
	$\operatorname{conv}c_H$. This proves the geometric hypothesis
	of Theorem~\ref{DFE} with $\omega=\lambda p$ and $\kappa=y_p$
	for every $p\in q(K)$. Its upper semicontinuity hypothesis follows
	from the continuity of $q$.
	
	Theorem~\ref{DFE} yields $\bar x\in K$ such that
	\[
	\langle u-\bar x,\nabla\Phi_{q(\bar x)}(u)\rangle\ge0
	\qquad(u\in\mathbb R^d).
	\]
	For $h\ne0$, the function
	$t\mapsto\Phi_{q(\bar x)}(\bar x+th)$ is therefore nondecreasing
	on $[0,\infty)$. It is real analytic and tends to $+\infty$ as
	$t\to+\infty$, so it cannot be constant on any nondegenerate
	interval. This proves the assertion.
\end{proof}

\begin{remark}
	Take $K=[0,\sqrt3]$, $\sigma=1$, $(b_1,b_2)=(0,\sqrt3)$,
	$(w_1,w_2)=(9,1)$, $\lambda=1/8$, and $q(x)=\sqrt3-x$.
	Then $D(H)=\mathbb R$ and
	\[
	H'(\sqrt3)=\tfrac18+9(-\tfrac14)+2=-\tfrac18.
	\]
	By continuity, $H'<0$ at interior points of $K$ sufficiently close
	to $\sqrt3$. Thus even $H|_K$ is nonmonotone, and neither
	Theorem~\ref{BR} nor Theorem~\ref{DF} applies directly.
\end{remark}
\eject

\end{document}